\documentclass[10pt,oneside]{article}
\usepackage{times,amsmath,amsbsy,amssymb,amscd,mathrsfs,graphicx}
\usepackage{amssymb}
\usepackage{amsmath,bm}
\usepackage{amsmath,amsthm}
\usepackage{pifont}
\usepackage{subcaption}
\usepackage{framed}
\usepackage{xcolor,color}
\usepackage[font=small,labelfont=md,textfont=it]{caption}
\usepackage{floatrow,float}
\usepackage[titletoc, title]{appendix}
\usepackage[colorlinks,linkcolor=blue,citecolor=blue]{hyperref}
\usepackage{empheq}
\usepackage{longtable}
\usepackage{pgfplots}
\usepackage{diagbox}
\usepackage{booktabs,makecell,multirow}
\usepackage[capitalise, nosort]{cleveref}
\usepackage{algorithm}
\usepackage{algorithmic}
\usepackage{amsmath,bm}
\usepackage{cases}
\usepackage{geometry}
\usepackage{extarrows}
\usepackage{tcolorbox}
\usepackage{bbm}
\usepackage{syntonly}
\usepackage[figuresright]{rotating}
\usepackage{epstopdf}

\crefname{equation}{}{}
\crefname{lem}{Lemma}{Lemmas}
\crefname{tab}{Table}{Tables}
\crefname{fig}{Figure}{Figures}
\crefname{thm}{Theorem}{Theorems}
\crefname{assum}{Assumption}{Assumptions}
\crefname{eg}{Example}{Examples}

\newcommand{\proj}[0]{ {\bf proj}}

\newcommand{\dd}{\,{\rm d}}

\newcommand{\R}{{\mathbb R}}

\newcommand{\dual}[1]{\left\langle {#1} \right\rangle}

\newcommand{\argmin}[0]{ {\mathop{{\rm  argmin}}}}

\newcommand{\nm}[1]{\left\lVert {#1} \right\rVert}
\newcommand{\snm}[1]{\left\lvert {#1} \right\rvert}

\newcommand{\ssnm}[1]
{
	\left\vert\kern-0.25ex
	\left\vert\kern-0.25ex
	\left\vert
	{#1}
	\right\vert\kern-0.25ex
	\right\vert\kern-0.25ex
	\right\vert
}

\makeatletter
\def\spher@harm#1{%
	\vbox{\hbox{%
			\offinterlineskip
			\valign{&\hb@xt@2\p@{\hss$##$\hss}\vskip.2ex\cr#1\crcr}%
		}\vskip-.36ex}%
}
\def\gshone{\spher@harm{.}}
\def\gshtwo{\spher@harm{.&.}}
\def\gshthree{\spher@harm{.&.&.}}
\let\gsh\spher@harm
\makeatother

\usepackage{tabularx}
\newtheorem{coro}{Corollary}[section]

\newtheorem{assum}{Assumption}

\newtheorem{lem}{Lemma}[section]

\newtheorem{rem}{Remark}[section]

\newtheorem{thm}{Theorem}[section]
\newtheorem{theorem}{Theorem}[section]

\newtheorem{eg}{Example}

\newcolumntype{I}{!{\vrule width 1,5pt}}
\newlength\savedwidth

\newlength\savewidth

\newcounter{mnote}
\let\oldmarginpar\marginpar
\renewcommand\marginpar[1]
{\-\oldmarginpar[\raggedleft\footnotesize #1]{\raggedright\footnotesize #1}}

\makeatletter\def\@captype{table}\makeatother

\newfloatcommand{capbtabbox}{table}[][\FBwidth]
\usepackage[T1]{fontenc}
\usepackage[utf8]{inputenc}
\usepackage{authblk}
\usepackage[all]{xy}

\newcommand{\conv}[1]{{\bf conv}\left\{ {#1} \right\}}

\usepackage{graphicx}
\usepackage{tikz}
\tikzset{elegant/.style={smooth,thick,samples=50,black}}
\tikzset{eaxis/.style={->,>=stealth}}

\usepackage[misc]{ifsym}

\providecommand{\HAMGRevisionColor}{black}

\begin{document}
	\title{\Large \bf  An accelerated multiobjective gradient flow with implicit Hessian-driven damping
	\thanks{This work was supported by the National Natural Science Foundation of China (Grant No. 12401402), the Science and Technology Research Program of Chongqing Municipal Education Commission (Grant No. KJZD-K202300505), the Natural Science Foundation of Chongqing (Grant No. CSTB2024NSCQ-MSX0329) and the Foundation of Chongqing Normal University (Grant No. 22xwB020).}}

\author[,1]{Wen Jiang\thanks{Email: 2024110518030@stu.cqnu.edu.cn}}
\author[,1,2]{Hao Luo \thanks{Email: luohao@cqnu.edu.cn}}
\affil[1]{National Center for Applied Mathematics in Chongqing, Chongqing Normal University, Chongqing, 401331, China}
\affil[2]{Chongqing Research Institute of Big Data, Peking University,  Chongqing, 401121, China}

	\date{\today}
	\maketitle
	
		\begin{abstract}
	In this paper, we propose an accelerated multiobjective gradient flow with implicit Hessian-driven damping for unconstrained convex multiobjective optimization problems in a real Hilbert space. The continuous flow is formulated as a system of two first-order equations involving only a single multiobjective projection operator and suitable time-scaling factors. Formally, it can also be expressed as a second-order system with explicit multiobjective Hessian-driven damping, which corresponds to the usual Hessian-driven damping introduced by Alvarez et al. [J. Math. Pures Appl., 81(8):747--779, 2002]. The first-order formulation, however, involves no Hessian term and can be regarded as an extension of the Hessian-driven Nesterov accelerated gradient flow by Chen and Luo [arXiv:1912.09276, 2019]. Assuming that the objective gradients are Lipschitz continuous on bounded sets and that the common lower level sets are bounded, we establish global existence in finite dimensions under suitable time-scaling conditions. An a priori energy estimate confines the solutions to a fixed bounded set, and a Lyapunov argument yields fast decay of a merit function and the KKT residual. Under an additional bound on the time-scaling factors, the trajectory of every global solution converges weakly to a weak Pareto point. Numerical tests illustrate the rapid residual decay and reduced oscillations.
\end{abstract}

		
\section{Introduction}
Let $\mathcal H$ be a real Hilbert space, possibly infinite-dimensional, with inner product $\dual{\cdot,\cdot}$ and induced norm $\nm{\cdot}$. In this paper, we consider the unconstrained multiobjective optimization problem
\begin{equation}\label{eq:min-mobj}
	\tag{MOP}
	\min_{x\in\mathcal H} F(x)
	=
	\left(f_1(x),\ldots,f_m(x)\right)^\top,
\end{equation}
where for $1\leq j\leq m$, each \(f_j:\mathcal H\to\R\) is convex and continuously differentiable, and the precise smoothness assumptions are stated in \cref{assum:f_j}.
Unlike the scalar case, the objective in \cref{eq:min-mobj} is now a vector-valued function with possibly conflicting components. Hence one cannot expect generally the existence of a point minimizing all objective components simultaneously. The appropriate solution concept is therefore given by Pareto or weak Pareto optimality, which makes both the theoretical analysis and numerical treatment more involved \cite{Ehrgott2005}. Classical methods typically use scalarization or heuristics. Scalarization reduces \cref{eq:min-mobj} to a single objective optimization problem, but its performance depends on the choice of scalar parameters. Heuristic methods can explore the Pareto front more flexibly, but often lack rigorous convergence guarantees.

Multiobjective gradient methods \cite{Drummond2004,fliege_steepest_2000,tanabe2019proximal,Tanabe2023a} provide a way to decrease all objective components simultaneously. A key ingredient is the multiobjective steepest descent direction \cite{fliege_steepest_2000}, which corresponds to the negative minimum-norm element $d = -\proj_{C(x)}(0)$ of the convex hull $C(x):=
\conv{\nabla f_1(x),\cdots,\nabla f_m(x)}$.  
From the continuous perspective, the dynamical system approach \cite{attouch2015multiibjective,attouch_multiibjective_2015,attouch2014continuous} gives an alternative framework, which focuses on continuous gradient dynamics driven by the multiobjective steepest descent direction. Thus, both discrete and continuous approaches involve projection onto the convex hull of gradients.
\subsection{Related works}
\subsubsection{Dynamical approach for the single objective case}
The dynamical system approach offers a continuous-time framework for unconstrained single objective optimization problems and has important connections with first-order methods. Beginning with Polyak's heavy-ball dynamics~\cite{Polyak1964}, this line of research was further developed by Su et al.~\cite{su_dierential_2016}, who derived a continuous-time model (also known as the asymptotic vanishing damping model \cite{Attouch2018}) for Nesterov's accelerated gradient method~\cite{Nesterov1983}. Since then, the dynamical approach and Lyapunov analysis have become useful tools for first-order methods; see
\cite{attouch2019rate,chen2025accelerated,chen2021unified,li2024linear,Luo2022b,Luo2023,Luo2024ax,Luo2024b,luo_accelerated_2026,Luo2025y,
	luo2025unified,
	Shi2022,wibisono2016variational} and the references therein.

However, due to the inertial (or momentum) terms, both continuous and discrete trajectories exhibit dramatic oscillations. To overcome this issue,
Alvarez et al.~\cite{alvarez2002second} first introduced the so-called Hessian-driven
damping and proposed the dynamical inertial Newton system
\begin{equation}\label{eq:din}
	\tag{DIN}
	x''(t)
	+ \alpha x'(t)
	+\beta\nabla^{2}F\bigl(x(t)\bigr)\cdot x'(t)
	+\nabla F\bigl(x(t)\bigr)
	=0,
\end{equation}
where $\alpha,\,\beta>0$ are constants. The Hessian-driven
damping involves a geometric mechanism governed by the Hessian of the
objective function and is naturally related to Newton-type methods.  Attouch et
al.~\cite{attouch2016fast} combined Hessian-driven damping with
asymptotically vanishing damping and introduced the DIN-AVD
system:
\begin{equation}
	\tag{DIN-AVD}\label{eq:DIN-AVD}
	x''(t)
	+\frac{\alpha}{t}x'(t)
	+\beta\nabla^{2} F\bigl(x(t)\bigr)\cdot x'(t)
	+\nabla F\bigl(x(t)\bigr)
	=0.
\end{equation} 
Compared to the constant damping \cite{Polyak1964} and the asymptotic vanishing damping \cite{su_dierential_2016}, the additional
Hessian-driven damping can substantially attenuate oscillations; see concrete numerical evidence in
\cite[Figure~1]{alvarez2002second}, \cite[Figure 1]{attouch2022first} and \cite[Section 1.5.5]{attouch2016fast}. One theoretical indication of this effect is that for any smooth convex objective $F$ with a nonempty minimizer set, \cref{eq:DIN-AVD} with $\alpha\geq3$ and $\beta>0$ admits the weighted gradient estimate \cite[Proposition 1.11]{attouch2016fast}
\begin{equation}\label{eq:fast-kkt}
	\int_{t_0}^{\infty}t^2\nm{\nabla F(x(t))}^2\dd t  <+\infty,
\end{equation}
where $t_0>0$ is the initial time.

\begin{figure}[H]
	\centering
	\begin{minipage}{0.5\textwidth}
		\centering
		\includegraphics[width=\linewidth]{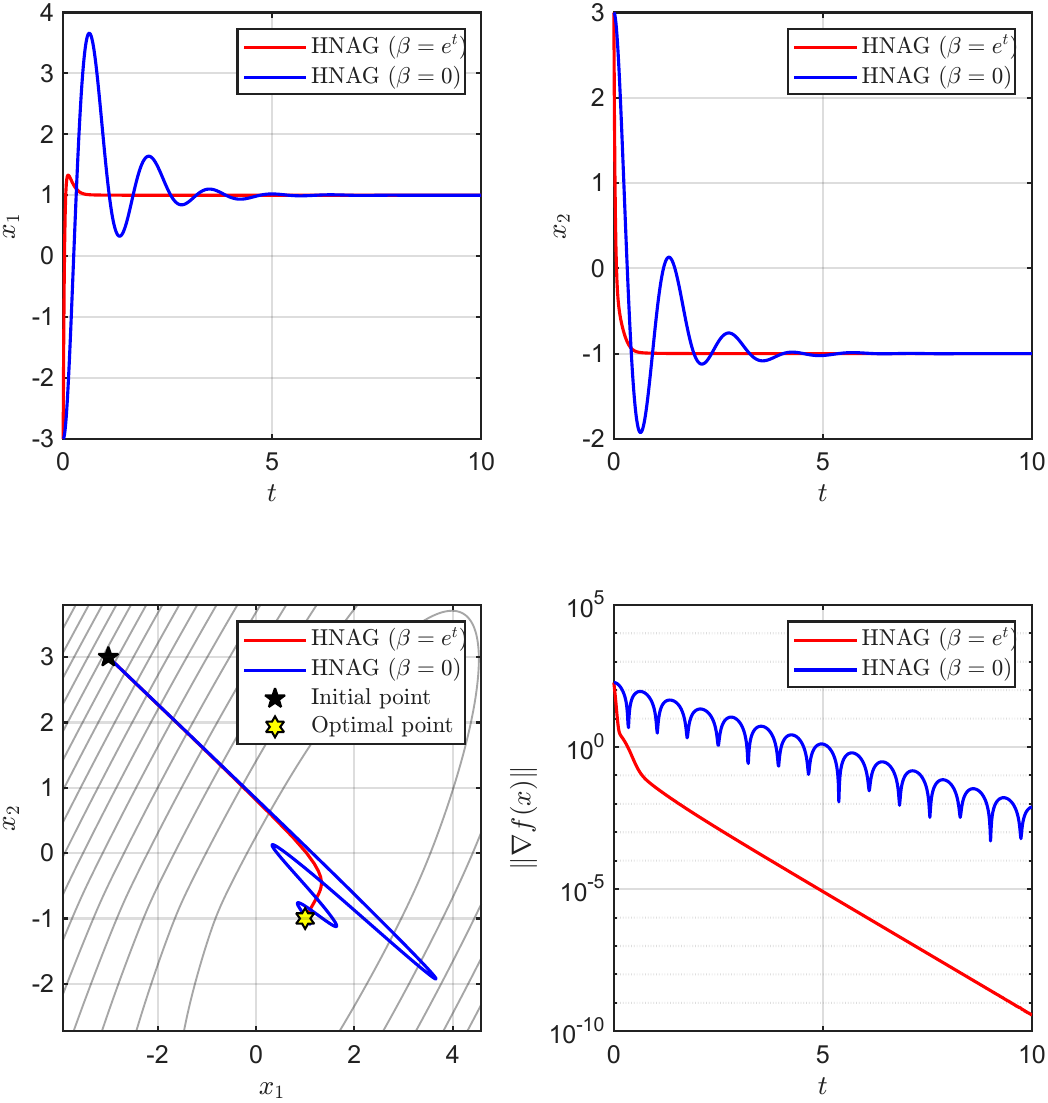}
	\end{minipage}
	\caption{Convergence behavior of \cref{eq:HNAG_second} for minimizing the log-sum-exp objective in \cref{eg:hessian}. Top line: the components of the trajectory. Bottom line: the trajectory in the \((x_1,x_2)\)-plane (left) and the KKT residual (right).}
	\label{fig:hessian}
\end{figure}

Building on the Nesterov accelerated gradient (NAG) flow~\cite{Luo2022d}, Chen and Luo~\cite{chen_first_2019} introduced the
Hessian-driven Nesterov accelerated gradient (HNAG)  flow
\begin{equation}\label{eq:HNAG_second}
	\tag{HNAG}
	\gamma(t) x''(t)
	+(\mu+\gamma(t))x'(t)
		+\gamma(t) \beta(t) \nabla^2 F(x(t))\cdot x'(t)
	+b(t)\nabla F(x(t))
	=0,
\end{equation}
where $b(t):=1+\mu\beta(t) +\gamma(t) \beta'(t)$, $\beta:[0,\infty)\to\R_+$ is a nonnegative $C^1$ function and the scaling factor $\gamma:[0,\infty)\to\R_+$ satisfies $	\gamma'(t)=\mu-\gamma(t)$,
with $\mu\geq0$ being the convexity parameter of the objective $F$. Note that \cref{eq:HNAG_second} can be reformulated as a first-order system without Hessian:
\begin{equation} \label{eq:HNAG-first-order}
	\left\{
	\begin{aligned}
	{}&	x'(t)
		 =
		v(t)-x(t)-\beta(t)\nabla F(x(t)),
		\\
	{}&	\gamma(t) v'(t)
		 =
		\mu(x(t)-v(t))-\nabla F(x(t)).
	\end{aligned}
	\right.
\end{equation}
For \(\beta=0\), this gives the NAG flow~\cite{Luo2022d} and the Hessian-driven damping disappears. For $\beta\neq0$, according to \cite[Theorem 2.1]{chen_first_2019}, the HNAG flow also satisfies a weighted gradient residual estimate like \cref{eq:fast-kkt}. As illustrated by \cref{eg:hessian}, the oscillation-reduction effect
can also be observed from the trajectory of \cref{eq:HNAG_second} with $\beta\neq0$.
\begin{eg}\label{eg:hessian}
Consider the log-sum-exp objective 
\[
	F(x)
=
\frac{\mu}{2}\|x\|^2
+
\ln
\left(
\sum_{i=1}^{p}
\exp
\left(
\langle a_i,x\rangle-b_i
\right)
\right),\quad x\in\R^2,
\]
where 
	\(a_i\in\mathbb{R}^2\) and \(b_i\in\mathbb{R}\).
We choose
	\(p=50\) and \(\mu=0.05\), and the entries of \(a_i\) and the scalars \(b_i\) are sampled from the uniform distribution on \([-1,1]\). Numerical results are reported in \cref{fig:hessian}.
\end{eg}

\subsubsection{Multiobjective dynamical models}

For \cref{eq:min-mobj}, the dynamical approach can be traced back
to the early works of Henry~\cite{Henry1973}
and Smale~\cite{Smale1973}. The underlying idea is mainly based on the multiobjective steepest descent direction, possibly combined with inertial and damping terms. Most existing multiobjective dynamical systems admit the following unified form
\begin{equation*}
	\theta X''(t) + \xi X'(t)
	+ \proj_{C(X(t))+Y}(\Delta) = 0,
\end{equation*}
where $\theta,\xi$ are parameters and $Y,\Delta$ are proper corrections. We refer to \cref{tab:multitarget_gradient} for a list of related works.

The multiobjective gradient system
\begin{equation}\label{eq:MOG}
	\tag{MOG}
	X'(t) + \proj_{C(X(t))}(0) = 0
\end{equation}
has been analyzed by Miglierina \cite{Miglierina2004}, Garrigos \cite{Garrigos2015}, and Attouch et al. \cite{attouch_multiibjective_2015,attouch2014continuous}. It can be regarded as the continuous-time limit of several
first-order multiobjective methods, including the steepest descent method \cite{drummond2005steepest,fliege_steepest_2000} and proximal gradient methods \cite{Bonnel2005,Drummond2004,tanabe2019proximal}. 
The inertial variant of \cref{eq:MOG} was proposed by Attouch and Garrigos \cite{attouch2015multiibjective}:
\begin{equation}\label{eq:IMOG}
	\tag{IMOG}
	\theta X''(t) + \xi X'(t) + \proj_{C(X(t))}(0) = 0,
\end{equation}
which  corresponds to the heavy-ball dynamics~\cite{Polyak1964}. The convergence analysis, however, requires the trajectory to be bounded and the constant damping parameter to
satisfy $\xi>\sqrt{\theta L}$, where $L$ denotes the common gradient Lipschitz constant  of all objectives along the bounded  trajectory.

\begin{table}[H]
	\centering
	\caption{Multiobjective gradient systems and their established properties. We use the abbreviations: "B.F.B." = "bounded from below", "B.L.L.S." = "bounded lower level set", "L.C.G." = "Lipschitz continuous gradients", "Conv." = "convergence". In the column of "Setting", $\mathcal H$ is a real Hilbert space, possibly infinite-dimensional; in the column of "L.C.G.", "Bounded" means Lipschitz continuous gradients on bounded sets while "Global" stands for globally Lipschitz continuous gradients; in the column of "Existence", an asterisk $*$ means an additional finite-dimensional assumption, and a star $\star$ means some additional technical assumptions; in the column of "Conv.", 	a dagger $\dagger$ indicates that
		boundedness of the trajectory is assumed, and a double dagger $\ddagger$ corresponds to some extra technical assumptions. }
	\label{tab:multitarget_gradient}
	
	\footnotesize
	\renewcommand{\arraystretch}{1.35}
	\setlength{\tabcolsep}{1.5pt}
	
	\begin{tabularx}{\linewidth}{@{}
			>{\hsize=1.15\hsize\linewidth=\hsize\raggedright\arraybackslash}X
			>{\hsize=0.6\hsize\linewidth=\hsize\centering\arraybackslash}X
			>{\hsize=1.4\hsize\linewidth=\hsize\centering\arraybackslash}X
			>{\hsize=1.5\hsize\linewidth=\hsize\centering\arraybackslash}X
			>{\hsize=1.27\hsize\linewidth=\hsize\centering\arraybackslash}X
			>{\hsize=0.66\hsize\linewidth=\hsize\centering\arraybackslash}X
			>{\hsize=0.60\hsize\linewidth=\hsize\centering\arraybackslash}X
			>{\hsize=0.81\hsize\linewidth=\hsize\centering\arraybackslash}X			
			>{\hsize=0.85\hsize\linewidth=\hsize\centering\arraybackslash}X				
			>{\hsize=0.95\hsize\linewidth=\hsize\centering\arraybackslash}X
			>{\hsize=1.47\hsize\linewidth=\hsize\centering\arraybackslash}X
			>{\hsize=0.7\hsize\linewidth=\hsize\centering\arraybackslash}X
			@{}}
		\toprule
		Model & $\theta$ & $\xi$ & $Y$ & $\Delta$
		& Setting & B.F.B. &B.L.L.S.&L.C.G.& Existence & Conv. & Rate\\
		\midrule
		
		MOG \cite{attouch2014continuous}
		& $0$ & $1$ & $0$ & $0$
		& $\mathcal H$ & Yes&No&Bounded
		& Global$^*$
		& Weak$^\dagger$
		& --\\
		
		IMOG \cite{attouch2015multiibjective}
		& $\theta>0$ & $\xi>\sqrt{\theta L}$ & $0$ & $0$
		& $\mathcal H$ & No&No&Bounded
		& Global$^{*,\star}$
		& Weak$^\dagger$
		& --\\
		
		IMOG' \cite{sonntag2024fast}
		& $1$ & $\xi>0$ & $0$ & $-X''(t)$
		& $\mathcal H$ & No&No&Bounded
		& Global$^{*,\star}$
		& Weak$^\dagger$
		& --\\
		
		MAVD \cite{Sonntag2024a}
		& $1$
		&  $\frac{\alpha}{t},\alpha\geq 3$
		& $0$ & $-X''(t)$
		& $\mathcal H$ & Yes&No&Global
		& Global$^{*}$
		&  { Weak$^\dagger({\alpha>3})$}
		& $\mathcal O(\frac1{t^{2}})$\\
		
		MAVD \cite{2025Yin}
		& $1$ & $3/t$ & $0$ & $-X''(t)$
		& $\R^n$ & Yes&Yes&Global
		& Global
		& Strong
		& $\mathcal O(\frac1{t^{2}})$\\
		
		MTRIGS\cite{boct2025inertial}
		& $1$
		& \parbox[c]{\linewidth}{\centering
			$\frac{\alpha}{t^q}, \alpha>0$\\
			$0<q<1$}
		& \parbox[c]{\linewidth}{\centering
			${ \frac{\beta}{t^p}X(t), \beta>0}$\\
			${ 0<p<q+1}$}
		& $-X''(t)$
		& $\mathcal H$ & No&No&Global
		& Global$^*$
		& Strong$^{\dagger,\ddagger}$
		& $\mathcal O(\frac1{t^{p}})$\\
		
		AMG \cite{Luo2025a}
		& $\gamma(t)$ & $\mu+\gamma(t)$ & $0$
		& $-\gamma(t)X''(t)$
		& $\R^n$ & No&Yes&Global
		& --
		& --
		& $\mathcal O(\frac1{e^{t}})$\\
		
		\bottomrule
	\end{tabularx}
\end{table}

Recently, Sonntag and Peitz~\cite{sonntag2024fast} introduced the inertial gradient-like system
\begin{equation}\label{eq:IMOG'}
	\tag{IMOG'}
	X''(t)+\xi X'(t)+\proj_{C(X(t))}\bigl(-X''(t)\bigr)=0,
\end{equation}
whose bounded trajectories converge weakly for any constant damping $\xi>0$. Sonntag and Peitz~\cite{Sonntag2024a} also introduced
a multiobjective inertial gradient-like system with asymptotically
vanishing damping:
\begin{equation}\label{eq:MAVD}
	\tag{MAVD}
	X''(t)+\frac{\alpha}{t} X'(t)
	+\proj_{C(X(t))}\bigl(-X''(t)\bigr)=0,
\end{equation}
which has the decay rate $\mathcal O(1/t^2)$ for $\alpha\geq3$ and weak convergence of bounded trajectories for $\alpha>3$. For $\alpha=3$, Yin \cite{2025Yin} established point convergence in $\R^n$. Using vanishing Tikhonov regularization and slowly vanishing damping, Bo\c{t} and Sonntag~\cite{boct2025inertial} proved strong convergence of bounded trajectories.
Compared with \cref{eq:MOG,eq:IMOG}, both \cref{eq:IMOG',eq:MAVD} adopt the inertial  term \(-X''(t)\) as a shifted projection argument.

Luo et al.~\cite{Luo2025a}
showed that  the continuous-time limit of the multiobjective accelerated
proximal gradient method by Tanabe et al.~\cite{Tanabe2023a} agrees with \cref{eq:MAVD}. They introduced an accelerated multiobjective
gradient flow with adaptive time scaling:
\begin{equation}\label{eq:AMG}
	\tag{AMG}
	\gamma(t)X''(t)+\bigl(\mu+\gamma(t)\bigr)X'(t)
	+\proj_{C(X(t))}\bigl(-\gamma(t)X''(t)\bigr)=0,
\end{equation}
where $\gamma'(t)=\mu-\gamma(t)$ with $\gamma(0)=\gamma_0>0$ and  $\mu\geq0$ is a common strong-convexity parameter. This provides a unified manner for convex ($\mu=0$) and strongly convex ($\mu>0$) problems. We refer to \cite{Luo2025v} for a recent primal-dual extension of \cref{eq:AMG} to linearly constrained multiobjective optimization problems.
\subsection{Motivation}
Among existing multiobjective gradient methods and dynamical systems, we mention some limitations.
\begin{itemize}
	\item[\textbf{(i)}] \textbf{Projection with inertial term.}
	Instead of the multiobjective steepest descent direction \cite{attouch2015multiibjective,attouch_multiibjective_2015,attouch2014continuous}, recent works 			\cite{boct2025inertial,Luo2025a,Sonntag2024a,sonntag2024fast} adopt projections of some special corrections depending on the inertial term (second-order derivative), which makes the theoretical
	analysis more involved, as the existence proofs are mainly based on equivalent differential inclusions.
\item[\textbf{(ii)}] \textbf{Fast decay on the KKT residual.}
Most existing works derive only the convergence rates with respect to merit functions \cite{Tanabe2024}, which, however, are not computable in practice. Although the Karush--Kuhn--Tucker (KKT) residual (cf.\cref{eq:1st-opt-cond-conv}) can be controlled by proper merit functions (cf.\cref{eq:u0-kkt}), it is rare to see a fast decay rate for the KKT residual itself like \cref{eq:fast-kkt}.
	
	\item[\textbf{(iii)}] \textbf{Inertia-induced oscillations.}
	The inertial (or momentum) terms introduced in multiobjective gradient methods and dynamical systems usually produce significant oscillations, as demonstrated by the numerical results in \cite{boct2025inertial,Luo2025a,Sonntag2024a,Tanabe2023a}. This breaks the monotonicity of the multiobjective gradient descent and deteriorates the practical performance. 
\end{itemize}
%
\subsection{Main contribution}
%
%

%

In this work, we propose a novel implicit Hessian-driven accelerated multiobjective
gradient (HAMG) flow
\begin{equation}\label{eq:HAMG-intro}
		\small
	\tag{HAMG}
	\left\{
	\begin{aligned}
		X'(t)
		={}&
		\eta(t)\sigma(t)(V(t)-X(t))
		-
		\eta(t)\sigma^{2}(t)\proj_{C(X(t))}\left(\frac{V(t)-X(t)}{\sigma(t)}\right),\\
		V'(t)
		={}	&
		-\mu\sigma(t)(V(t)-X(t))
		-
		\sigma(t)\proj_{C(X(t))}\left(\frac{V(t)-X(t)}{\sigma(t)}\right),
	\end{aligned}
	\right.
\end{equation}
where 
$\eta'(t)
=
\sigma(t)		\eta(t)\bigl(\mu-\eta(t)\bigr)$ and
\(\sigma:[0,\infty)\to(0,\infty)\) is a continuous function, with the initial conditions $X(0)=X_0\in\mathcal H,\,V(0)=V_0\in\mathcal H$ and $\eta(0)=\eta_0>0$.  Again, $\mu\geq0$ is a common strong-convexity parameter of all objective components. Unlike most existing multiobjective dynamical systems in \cref{tab:multitarget_gradient}, our \cref{eq:HAMG-intro} does not involve the projection of the inertial \(X''(t)\). If $\sigma\in C^1$, \cref{eq:HAMG-intro} can formally be rewritten as 
\begin{equation}\label{eq:HAMG-X-2nd-intro}	
	X''(t) + \frac{\alpha(t)}{\sigma(t)}X'(t)
	+ \eta(t)\sigma^2(t)	\frac{\dd}{\dd t}
	\proj_{C(X(t))}(U(t))
	+ \rho(t)		\proj_{C(X(t))}(U(t)) = 0,
\end{equation}
where $\alpha(t)=2\eta(t)\sigma^2(t) - \sigma'(t) ,\, \rho(t)=\eta(t)\sigma(t)\left(\sigma(t)+\mu\sigma^2(t) + \sigma'(t) \right)$ and $U(t) = (V(t)-X(t))/\sigma(t)$. As one can see from \cref{eq:HAMG-X-2nd-intro}, there is an explicit differential term of the projection, which corresponds to the Hessian-driven damping in \cref{eq:din,eq:HNAG_second} for single-objective optimization problems. Thus, we might expect that  \cref{eq:HAMG-intro} exhibits less oscillatory trajectories and faster residual decay, as illustrated numerically by \cref{eg:small_pareto_set}. This is consistent with the oscillation-reduction effect of the Hessian-driven damping in \cref{eq:din,eq:HNAG_second}.

\begin{figure}[h]
	\centering
	\begin{minipage}{0.5\textwidth}
		\centering
		\includegraphics[width=\linewidth]{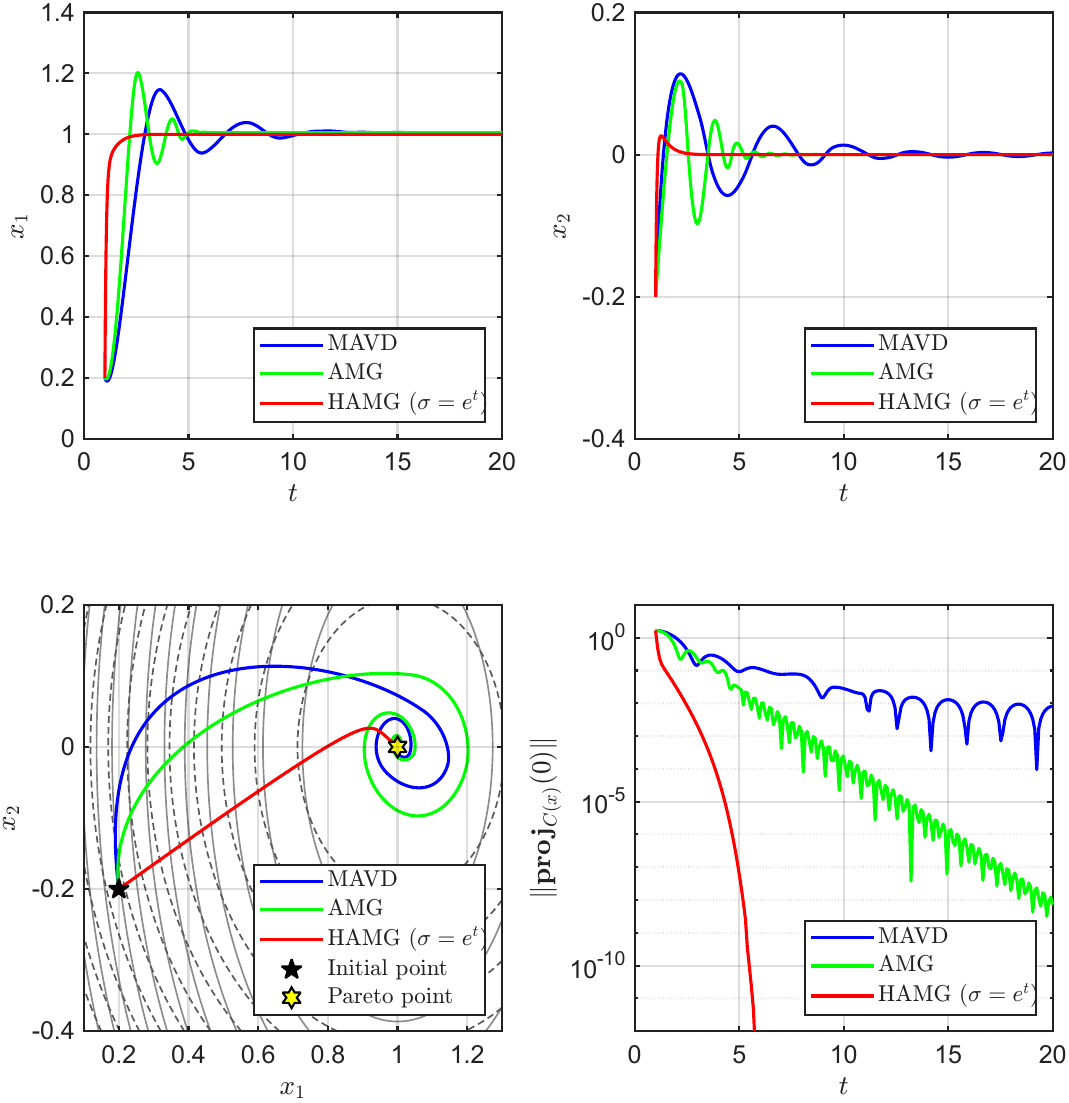}
	\end{minipage}
	\caption{Convergence behaviors of \cref{eq:HAMG-intro}, \cref{eq:MAVD,eq:AMG} for \cref{eg:small_pareto_set}. Top line: the components of the trajectory. Bottom line: the trajectory in the \((x_1,x_2)\)-plane (left) and the KKT residual (right).}
	\label{fig:hessian_mop}
\end{figure}

\begin{eg}\label{eg:small_pareto_set} 
	Consider a two-dimensional bi-objective optimization problem 
	\[
	\begin{aligned}
		f_1(x)
		&=
		(x_1-1)^2+x_2^2, \quad
		f_2(x)
		=
		2(x_1-1.01)^2+4x_2^2.
	\end{aligned}
	\]
	The (weak) Pareto set is 
	$ 
	\mathcal P_{w} 
	= 
	\left\{ 
	(1+0.01\theta,\,0)^\top: 
	0\leq \theta\leq1
	\right\}
	$ and a common strong-convexity parameter is $\mu=2$. The numerical results have been reported in  \cref{fig:hessian_mop}.
\end{eg}

Under \cref{assum:f_j,assum:L,assum:eta-sigma}, \cref{eq:HAMG-intro} has a global $C^1$ solution $(X,V):[0,\infty)\to\mathcal H\times\mathcal H$ when $\dim\mathcal H<\infty$. 
For any given global $C^1$ solution, we derive the following decay estimate
\begin{equation}\label{eq:decay-intro}
	\begin{aligned}
		u_0(X(t))+ \frac1{\Sigma(t)} \int_{0}^{t}\sigma^{2}(\tau)\eta(\tau)\Sigma(\tau)\big\|\proj_{C(X(\tau))}(0)\big\|^2 \dd \tau\leq
		\mathcal O\left(	\frac{1}{\Sigma(t)}\right),
	\end{aligned}
\end{equation}
where $u_0$ is a merit function (cf.\cref{eq:merit}) and $\Sigma(t)=e^{\int_0^t\sigma(\tau)\eta(\tau)\dd\tau}\to\infty$ as $t\to\infty$. We refer to \cref{tab:beta_convergence_rate-} for the explicit growth behaviors of $\Sigma$ with specific choices of $\sigma$. In particular, we obtain the fast convergence rate on the KKT residual
\begin{equation}\label{eq:rate-kkt}
	\small
	\begin{aligned}
		\big\|\proj_{C(X(t))}(0)\big\|^2 \leq 	\mathcal O\left(	\frac{1}{\Sigma(t)}\right),\quad 	\inf_{0<\tau<t}
		\big\|\proj_{C(X(\tau))}(0)\big\|^2
		=
		\mathcal O\left(
		\frac{1}{1+
			\int_{0}^{t}
			\sigma^{2}(\tau)\eta(\tau)\Sigma(\tau)\dd \tau
		}
		\right).
	\end{aligned}
\end{equation}
Under extra restriction $\eta^2(t)\sigma(t)\leq2\mu+6\eta(t)$, we prove weak convergence of $X(t)$ and $V(t)$ to a weak Pareto optimal solution of \cref{eq:min-mobj}, and strong convergence if $\mu>0$.

We mention that \cref{eq:HAMG-intro} is mainly motivated by \cref{eq:HNAG-first-order}, i.e., the Hessian-free presentation of \cref{eq:HNAG_second}. 
As discussed in \cref{sec:rescal}, when $\sigma(t)\eta(t)=1$, 
\cref{eq:HAMG-intro} gives a new multiobjective gradient system  with asymptotic vanishing anchor
\begin{equation}  \label{eq:hamg-mog-ava-intro}
	\tag{MOG-AVA}
	\eta(t)X'(t) =	\eta_0e^{-2t}(V_0-X_0)-  \proj_{C(X(t))}(\eta_0e^{-2t}(V_0-X_0)),
\end{equation} 	
for which the decay estimate \cref{eq:decay-intro} holds with $\Sigma(t) = e^t$. 
We show further that, when $\eta\sigma=1$, under the rescaling $t(s)=\int_{0}^{s}\lambda(r)\dd r$ with $s\geq0$ (cf.\cref{eq:res-2-din}),  	\cref{eq:HAMG-intro} becomes
\begin{equation}
	\label{eq:hamg-x-v-intro}
	\left\{
	\begin{aligned}					\lambda(s)x'(s)={}&\lambda^2(s)\big(v(s)-x(s)\big)-\mathbf{proj}_{C(x(s))}(u(s)),\\ \lambda(s)v'(s)={}&-\mu\bigl(v(s)-x(s)\bigr)-\mathbf{proj}_{C(x(s))}(u(s)),
	\end{aligned}
	\right.
\end{equation}
where $u(s):=\eta_0(V_0-X_0)e^{-2 \int_{0}^{s}\lambda(r)\dd r}$. 
Under the stated assumptions, the rescaled system \cref{eq:hamg-x-v-intro} admits a global $C^1$ solution in finite dimensions. For any global solution $(x,v):[0,\infty)\to\mathcal H\times\mathcal H$, both $x(s)$ and $v(s)$ converge weakly to a weak Pareto point of \cref{eq:min-mobj}, and the convergence is strong if $\mu>0$. 
Moreover, the second-order model \cref{eq:HAMG-X-2nd-intro} is transformed into
\begin{equation}\label{eq:hamg-x}
		\small
	\begin{aligned}
		x''(s) + \frac{\mu+\lambda^2(s)-\lambda'(s)}{\lambda(s)}x'(s)+\frac{1}{\lambda(s)}\frac{\dd}{\dd s}
		\proj_{C(x(s))}(u(s)) +2\proj_{C(x(s))}(u(s)) = 0.
	\end{aligned}
\end{equation}
This yields possible extensions of \cref{eq:din,eq:DIN-AVD} to the multiobjective setting. Specifically, for $\mu>0$ and $\eta_0=\mu$, we obtain from \cref{eq:hamg-x} a
multiobjective dynamical inertial Newton system
\begin{equation}\label{eq:hamh-din-intro}
	\tag{MDIN}
	x''(s) + 2\sqrt{\mu}x'(s)+\frac{1}{\sqrt{\mu}}\frac{\dd}{\dd s}
	\proj_{C(x(s))}(u(s)) +2\proj_{C(x(s))}(u(s)) = 0,
\end{equation}
with $u(s) = \mu (V_0-X_0)e^{-2\sqrt{\mu}s}$. Correspondingly, the decay estimates \cref{eq:decay-intro,eq:rate-kkt} become
\[
\begin{aligned}
	u_0(x(s))+ \frac1{\sqrt{\mu}} \int_{0}^{s}e^{\sqrt{\mu}(r-s)}\big\|\proj_{C(x(r))}(0)\big\|^2 \dd r\leq
	\mathcal O
	\left(e^{-\sqrt{\mu}s}\right),
\end{aligned}
\]
and 
\[
	\small
\begin{aligned}
	\big\|\proj_{C(x(s))}(0)\big\|^2 \leq 	\mathcal O\left(	e^{-\sqrt{\mu}s}\right),\quad 	\inf_{0<r<s}
	\big\|\proj_{C(x(r))}(0)\big\|^2
	=
	\mathcal O\left(
e^{-\sqrt{\mu}s}
	\right).
\end{aligned}
\]
For $\mu=0$ and $\eta_0=4$, we get a multiobjective dynamical inertial Newton system with asymptotic vanishing damping
\begin{equation}\label{eq:hamh-din-avd-intro}
		\small
	\tag{MDIN-AVD}
	x''(s) + \frac{3x'(s)}{s+1}+\frac{s+1}{2}\frac{\dd}{\dd s}
	\proj_{C(x(s))}(u(s)) +2\proj_{C(x(s))}(u(s)) = 0,
\end{equation}
with $u(s) = 4(V_0-X_0)(s+1)^{-4}$. In this case, the decay estimates \cref{eq:decay-intro,eq:rate-kkt} become
\[
\begin{aligned}
	u_0(x(s))+  \int_{0}^{s}\frac{(r+1)^3}{2(s+1)^{2}}\big\|\proj_{C(x(r))}(0)\big\|^2 \dd r\leq
	\mathcal O\left(	(s+1)^{-2}\right),
\end{aligned}
\] 
and 
\[
\small
\begin{aligned}
	\big\|\proj_{C(x(s))}(0)\big\|^2 \leq 	\mathcal O\left(	(1+s)^{-2}\right),\quad 	\inf_{0<r<s}
	\big\|\proj_{C(x(r))}(0)\big\|^2
	=
	\mathcal O\left(
	(1+s)^{-4}
	\right).
\end{aligned}
\]

\subsection{Organization} 

The remainder of this paper is organized as follows. We present some preliminaries in \cref{sec:pre} and in \cref{sec:hamg}, we focus on the dynamical system \cref{eq:HAMG-intro} and establish a basic energy estimate. Then in \cref{sec:exist}, we derive the local and global existence results in a finite-dimensional setting. In \cref{sec:fast}, we  develop a Lyapunov analysis and derive the convergence rate estimates of a merit function and the KKT residual. After that, we establish weak convergence of the trajectory and the corresponding strong convergence results in \cref{sec:x}. Finally, some numerical experiments are reported in \cref{sec:tests}, and concluding remarks are given in \cref{sec:conclu}.  

\section{Preliminaries}
\label{sec:pre}

This section collects some notation and preliminary results for later use. 
\subsection{Functional classes}
Denote the Hilbert inner product by \(\dual{\cdot,\cdot}\) and the induced norm by \(\nm{\cdot}\). 
For a function \(f:\mathcal H\to\R\),
define its lower level set by
$
\mathcal L_f(\alpha)
:=
\{x\in\mathcal H:\ f(x)\le \alpha\}
$, where  \(\alpha\in\R\).
For a vector-valued function
\(F=(f_1,\ldots,f_m)^\top:\mathcal H\to\R^m\), let
$
\nabla F(x)
:=
[\nabla f_1(x),\ldots,\nabla f_m(x)]$
be the adjoint of the Jacobian of $F$. For \(a\in\R^m\), define the common lower level set by
$
\mathcal L_F(a)
:=
\{x\in\mathcal H:f_j(x)\leq a_j,\,1\leq j\leq m\} = \cap_{1\leq j\leq m}\mathcal L_{f_j}(a_j).
$

For any $R>0$, write $\overline B_R:=\{x\in\mathcal H:\|x\|\leq R\}$.
Let $\mathcal F_{\mathrm{ball}}^{1,1}(\mathcal H)$ be the class of all continuously differentiable functions whose gradients are Lipschitz continuous on every closed ball (or equivalently bounded set): for each $R>0$, there is a finite constant $L_f(R)>0$ such that
\[
\|\nabla f(x)-\nabla f(y)\|\leq L_f(R)\|x-y\|,
\qquad \forall\,x,y\in\overline B_R.
\]
Let \(\mathcal S_\mu^1(\mathcal H)\) be the class of continuously differentiable \(\mu\)-strongly convex functions, where \(\mu\geq0\), namely,
\[
f(y)\geq f(x)+\dual{\nabla f(x),y-x}
+\frac{\mu}{2}\nm{y-x}^2,
\quad \forall\,x,y\in\mathcal H,
\]
which also yields, for every $\theta\in[0,1]$, the basic inequality
\begin{equation}\label{eq:scv}
	f(\theta x+(1-\theta)y)\leq \theta f(x)+(1-\theta)f(y)-\frac{\mu}{2}\theta(1-\theta)\nm{x-y}^2,\quad \forall\,x,y\in\mathcal H.
\end{equation}
It is standard that any $f\in \mathcal S_\mu^1(\mathcal H)$ is
 weakly lower semicontinuous.
For the subsequent analysis, we further set $\mathcal S_{\mu,\mathrm{ball}}^{1,1}(\mathcal H)
:=\mathcal F_{\mathrm{ball}}^{1,1}(\mathcal H)\cap\mathcal S_\mu^1(\mathcal H)$, and  impose the following two assumptions.
\begin{assum}\label{assum:f_j}
For every $1\leq j\leq m$, $f_j\in\mathcal S_{\mu_j,\mathrm{ball}}^{1,1}(\mathcal H)$ with $\mu_j\geq0$. We set $\mu:=\min_{1\leq j\leq m}\mu_j$.
\end{assum}
\begin{assum}\label{assum:L}
	The common lower level set $\mathcal L_F(a)$ is bounded for every \(a\in\R^m\).
	Equivalently, for any \(a\in\R^m\), $	S(a)
	:=
	\sup_{x\in\mathcal L_{F}(a)}
	\|x\|
	<+\infty$, whenever \(\mathcal L_{F}(a)\neq\emptyset\).
\end{assum}
\subsection{Weak Pareto optimality}

For \(a,b\in\R^m\), we write \(a\prec b\) if
\(a_j<b_j\) for all $1\leq j\leq m$ and 
\(a\preceq b\) if
\(a_j\le b_j\) for all $1\leq j\leq m$.
A point \(y\in\mathcal H\) is said to be dominated by
\(x\in\mathcal H\) with respect to \cref{eq:min-mobj} if $F(x)\preceq F(y)
$ and $
F(x)\ne F(y)$. Equivalently, we say that \(x\) dominates \(y\). A point \(x^\ast\in\mathcal H\) is called a Pareto optimal
solution or a Pareto point of \cref{eq:min-mobj} if there is no
\(y\in\mathcal H\) that dominates  $x^\ast$.
The set of all Pareto points of \cref{eq:min-mobj} is denoted by $\mathcal P
:=
\left\{
x^\ast\in\mathcal H:
\text{\(x^\ast\) is a Pareto point of \cref{eq:min-mobj}}
\right\}$. Its image \(F(\mathcal P)\) is called the Pareto front.  A point \(x^\ast\in\mathcal H\) is called a weakly Pareto optimal
solution or a weak Pareto point of \cref{eq:min-mobj} if there is no
\(y\in\mathcal H\) such that $F(y)\prec F(x^\ast)$.
The set of all weak Pareto points of \cref{eq:min-mobj} is denoted by $\mathcal P_w
:=
\left\{
x^\ast\in\mathcal H:
\text{\(x^\ast\) is a weak Pareto point of \cref{eq:min-mobj}}
\right\}$.
The image \(F(\mathcal P_w)\) is called the weak Pareto front. It is clear that $\mathcal P\subset\mathcal P_w$.
\begin{lem}[\cite{Sonntag2024a}]
	\label{lem:S}
	Let \(a\in\R^m\) be such that
	\(\mathcal L_F(a)\neq\emptyset\). Under \cref{assum:f_j,assum:L}, the common lower level set \(\mathcal L_F(a)\) is closed and convex, hence weakly closed and weakly compact. Moreover, we have $\mathcal P\cap\mathcal L_F(a)\neq\emptyset$ and $
	\sup_{b\in F(\mathcal P_w\cap\mathcal L_F(a))}
	\inf_{z\in F^{-1}(b)}
	\|z\|\leq S(a)
	<+\infty $.
\end{lem}

%

\subsection{Optimality condition}\label{subsec:kkt}
Define the convex hull of $\{p_1, \cdots, p_m\} \subset \mathcal H$ by $\conv{p_1, \cdots, p_m} := \{ \sum_{j=1}^m \lambda_j p_j : \lambda \in \Delta_m\}$, which is compact in $\mathcal H$. Here $\Delta_m := \{ \lambda \in \mathbb{R}^m_+ :  \sum_{j=1}^m \lambda_j=1\}$ denotes the unit simplex. 
A point \(x^\ast\in\mathcal H\) is said to satisfy the first-order
necessary optimality condition for \cref{eq:min-mobj} if
\begin{equation}\label{eq:1st-opt-cond-conv}
	\proj_{C(x^\ast)}(0)=0\quad\Longleftrightarrow\quad 	0\in C(x^\ast):=
	\conv{\nabla f_1(x^\ast),\cdots,\nabla f_m(x^\ast)},
\end{equation}
which is also referred to as the Karush--Kuhn--Tucker (KKT)
condition. Equivalently, there exists
$
\lambda^\ast\in\Delta_m
$
such that $	\nabla F(x^\ast)\lambda^\ast
=
\sum_{j=1}^m
\lambda_j^\ast\nabla f_j(x^\ast)
=
0$.
Thanks to \cref{assum:f_j}, the KKT
condition \cref{eq:1st-opt-cond-conv} is sufficient and necessary for weak Pareto optimality; see \cite{attouch_multiibjective_2015,Luo2025a,sonntag2024fast}.

\subsection{Merit function}
A merit function for \cref{eq:min-mobj} is a nonnegative function whose
zero set coincides with $\mathcal P_w$.
We adopt the following merit
function \cite{Liu2009}
\begin{equation}\label{eq:merit}
	u_0(x)
	:=
	\sup_{z\in\mathcal H}
	f(x;z),\quad \text{where}\,	f(x;z):=
	\min_{1\le j\le m}
	\left[f_j(x)-f_j(z)\right].
\end{equation}
According to \cite[Theorem~3.1]{Tanabe2024}, the function \(u_0:\mathcal H\to\R\cup\{+\infty\}\) defined
by \cref{eq:merit} is nonnegative and lower semicontinuous. Moreover, it provides a valid measure
of the weak Pareto optimality for \cref{eq:min-mobj}: 
\begin{equation}\label{eq:u_pw}
	x^\ast\in\mathcal P_w
	\quad\Longleftrightarrow\quad
	u_0(x^\ast)=0.
\end{equation}
%
For each $R>0$, choose the gradient Lipschitz constant $L_{f_j}(R)$ of the objective $f_j$ on $\overline B_R$ and define
\begin{equation}\label{eq:LR-MR}
	L_R:=\max\{L_{f_1}(R),\ldots,L_{f_m}(R)\},\qquad
	M_R:=\max_{1\leq j\leq m}\sup\{\|\nabla f_j(x)\|:\,x\in\overline B_R\}.
\end{equation}
These constants are finite, since $M_R\leq RL_R+\max_j\|\nabla f_j(0)\|$. 
Under \cref{assum:f_j}, we have the following local comparison, analogous to \cite[Theorem 4.1(ii)]{Tanabe2024}:
\begin{equation}\label{eq:u0-kkt}
u_0(x)\geq
\frac{1}{2N_R}\big\|\proj_{C(x)}(0)\big\|^2,
\quad \forall\,x\in\overline B_R,\quad R>0,
\end{equation}
with $N_R := \max\{L_{R+1},M_R\}$.
To see this,  let us fix any $x\in\overline B_R$ and set $p=\proj_{C(x)}(0)$. The projection optimality gives $\langle\nabla f_j(x),p\rangle\geq\|p\|^2$ for every $j$. Since $p\in C(x)$, we also have $\|p\|\leq M_R$. Thus $y=x-p/N_R$ belongs to $\overline B_{R+1}$. Applying the descent lemma on this larger ball yields
\[
f_j(y)\leq f_j(x)-\frac{\langle\nabla f_j(x),p\rangle}{N_R}
+\frac{L_{R+1}}{2N_R^2}\|p\|^2
\leq f_j(x)-\frac{\|p\|^2}{2N_R},\quad 1\leq j\leq m.
\]
Taking the minimum over $j$ and using $y$ as a competitor in \cref{eq:merit} proves the desired lower bound.

As a consequence of \cref{lem:S}, the supremum in \cref{eq:merit} can be restricted; see
 \cite[Lemma~2.2]{Luo2025a} and \cite[Theorem~6.7]{sonntag2024fast}. Additionally, this also yields that $u_0(x)<+\infty$ for all $x\in\mathcal H$.
\begin{lem}[ \cite{Luo2025a,sonntag2024fast}]\label{lem:L_F}
Let $a\in\R^m$ be such that $\mathcal L_F(a)\neq\emptyset$, then for every $x\in\mathcal L_F(a)$, we have $u_0(x)=\sup_{b\in F(\mathcal P_w\cap\mathcal L_F(a))}\inf_{z\in F^{-1}(b)}f(x;z)$.
\end{lem}


\subsection{Projection operator}
For any nonempty closed convex subset $C \subset \mathcal H$ and any point $x \in \mathcal H$, the distance from $x$ to $C$ is defined as $d(x, C) := \min_{z \in C} \| z - x \|$.
The minimizer exists uniquely and is denoted as $\proj_C(x) := \argmin_{z \in C}  \| z - x \|$, which satisfies 
\begin{equation}\label{eq:opt-proj}
	\dual{x-\proj_C(x),z-\proj_C(x)}\leq 0\quad\forall\,z\in C.
\end{equation}
It is easy to verify the following translation decomposition property:
\begin{equation*}
%
\proj_{x + C}(0) = x + \proj_C(-x) \quad \forall\, x \in \mathcal H.
\end{equation*}
In addition, $\proj_{C}(x)$ is 1-Lipschitz continuous with respect to $x$:
\begin{equation}\label{eq:Proj_Lipschitz}
\left\| \proj_{C} (x) - \proj_{C} (y) \right\|
\leq \nm{x - y} \quad \forall\, x ,y\in \mathcal H.
\end{equation}
Quantitative continuity estimates for projections under perturbations of the set are given in \cite[Example 4.1]{attouch1993quantitative}. 
%
The following local projection estimate is a bounded-ball counterpart of the estimates in \cite[Lemma 3.2]{Luo2025a} and \cite[Proposition 2.8]{attouch2014continuous}.
\begin{lem}\label{lem:proj_U}
Under \cref{assum:f_j}, for any $R>0$, $x,y\in\overline B_R$ and $z\in\mathcal H$, we have
\[
\big\|\proj_{C(x)}(z)-\proj_{C(y)}(z)\big\|^2
\leq2L_R\bigl(\|z\|+M_R\bigr)\|x-y\|.
\]
\end{lem}
\begin{proof}
Put $p=\proj_{C(x)}(z)=\nabla F(x)\lambda_x$ and $q=\proj_{C(y)}(z)= \nabla F(y)\lambda_y$, with $\lambda_x\in\Delta_m$ and $\lambda_y\in\Delta_m$.
Introduce  $\widehat p = \nabla F(x)\lambda_y\in C(x)$ and $\widehat q= \nabla F(y)\lambda_x\in C(y)$. Then it is clear that $\|\widehat p-q\|\leq\delta$ and $\|\widehat q-p\|\leq\delta$, with  $\delta=L_R\|x-y\|$. Applying \cref{eq:opt-proj} implies the following:
\[
\begin{aligned}
\|p-q\|^2
&=\langle z-p,q-p\rangle+\langle z-q,p-q\rangle
\leq\langle z-p,q-\widehat p\rangle+\langle z-q,p-\widehat q\rangle\\
&\leq\delta\bigl(\|z-p\|+\|z-q\|\bigr)
\leq2\delta(\|z\|+M_R),
\end{aligned}
\]
where the last step uses $\|p\|\leq M_R$ and $\|q\|\leq M_R$. This completes the proof.
\end{proof}

\section{HAMG: A New Multiobjective Gradient Dynamics}\label{sec:hamg}
\subsection{Continuous flow} 
Motivated by \cref{eq:HNAG-first-order}, we propose an implicit Hessian-driven accelerated multiobjective gradient (HAMG) flow
\begin{equation}\label{eq:HAMG}\small
	\tag{HAMG}
	\left\{
	\begin{aligned}
		X'(t)
		={}&
		\eta(t)\sigma(t)(V(t)-X(t))
		-
		\eta(t)\sigma^{2}(t)\proj_{C(X(t))}\left(\frac{V(t)-X(t)}{\sigma(t)}\right),\\
		V'(t)
		={}		&
		-\mu\sigma(t)(V(t)-X(t))
		-
		\sigma(t)\proj_{C(X(t))}\left(\frac{V(t)-X(t)}{\sigma(t)}\right),
	\end{aligned}
	\right.
\end{equation}
with the initial conditions $X(0) = X_0\in\mathcal H$ and $V(0) = V_0\in\mathcal H$. Here $\sigma:[0,\infty)\to\R_+$ is a positive and continuous
function and the scaling parameter $\eta$ satisfies
\begin{equation} 
	\eta'(t)
	={} 
	\sigma(t)		\eta(t)\bigl(\mu-\eta(t)\bigr),\quad \eta(0) = \eta_0>0.
	\label{eq:HAMG-eta}
\end{equation}
Recall that $\mu\geq0$ denotes a common strong-convexity parameter (cf.\cref{assum:f_j}).
The parameter equation \cref{eq:HAMG-eta} is a Bernoulli equation and it is not hard to obtain the explicit form
\begin{equation}\label{eq:eta}
	\eta(t) = \left\{
	\begin{aligned}
		& \frac{\eta_0}{1 + \eta_0 			\int_0^t \sigma(\tau)\dd \tau}, \quad && \mu = 0,\\
		& \frac{\mu\eta_0}{\eta_0 + (\mu-\eta_0)e^{-\mu 	\int_0^t \sigma(\tau)\dd \tau}}, \quad && \mu > 0.
	\end{aligned}
	\right.
\end{equation}
It is clear that $\min\{\eta_0,\mu\}=:\eta_{\min}\leq \eta(t)\leq \eta_{\max}:=\max\{\eta_0,\mu\}$ for all $t\geq0$.

For the following reformulations, assume additionally that $\sigma\in C^1$. Introduce the variable
\begin{equation*}
	U(t):=\frac{V(t)-X(t)}{\sigma(t)},\quad P(X(t),U(t)): = \proj_{C(X(t))}(U(t)),
\end{equation*}
then \cref{eq:HAMG} can be rewritten as 
\begin{subnumcases}{\label{eq:HAMG-X-U}} 
	X'(t)
	=
	\eta(t)\sigma^{2}(t)\bigl(U(t)-P(X(t),U(t))\bigr),
	\label{eq:HAMG-X-U-X}\\
	U'(t)
	=
	-\delta(t)U(t)
	+
	(\eta(t)\sigma(t)-1)P(X(t),U(t)),
	\label{eq:HAMG-X-U-U} 
\end{subnumcases}
with $\delta(t):=\sigma(t)(\mu+\eta(t))+\frac{\sigma'(t)}{\sigma(t)}$.
We obtain from \cref{eq:HAMG-X-U}  immediately that 
\begin{equation*}
	\eta(t)\sigma^2(t)U'(t) + \eta(t)\sigma(t)\left(\sigma(t)+\mu\sigma^2(t) + \sigma'(t) \right)U(t)
	= (1-\eta(t)\sigma(t))X'(t),
\end{equation*}
which yields the following relation
\begin{equation}\label{eq:U-def-second-order}
	U(t)
	=
	\frac{
		V_0-X_0
		+Z(t)}{\sigma(t)A(t)},\quad 	A(t)
	:	=
	\exp\left(
	\int_0^t
	\left(1+\mu\sigma(\tau)\right)\dd \tau
	\right),
\end{equation}
with
\begin{equation}\label{eq:Z-A}
	Z(t)
	:=
	\int_0^t\frac{A(\tau)}{\sigma(\tau)\eta(\tau)}
	\left(1-\sigma(\tau)\eta(\tau)\right)	 X'(\tau)\dd \tau.
\end{equation}

On the other hand, based on \cref{eq:HAMG-X-U}, eliminating $U$ and $U'$ gives formally a second-order model
\begin{equation}\label{eq:HAMG-X-2nd}	
	X''(t) + \frac{\alpha(t)}{\sigma(t)}X'(t)
	+ \eta(t)\sigma^2(t)	\frac{\dd}{\dd t}
	\proj_{C(X(t))}(U(t))
	+ \rho(t)		\proj_{C(X(t))}(U(t)) = 0,
\end{equation}
where $\alpha(t):=2\eta(t)\sigma^2(t) - \sigma'(t) $ and $ \rho(t):=\eta(t)\sigma(t)\left(\sigma(t)+\mu\sigma^2(t) + \sigma'(t) \right)$.
In \cref{eq:HAMG-X-2nd}, there is a special term involving the differentiation of the projection,
which is related to the well-known Hessian-driven damping  
in \cref{eq:DIN-AVD}  for the single objective case.  
However, the projection mapping
$
\proj_{C(X(t))}\bigl(U(t)\bigr)
$
may not be differentiable with respect to $t$. Thus, instead of the explicit representation in \cref{eq:HAMG-X-2nd}, the damping mechanism is understood implicitly through the first-order system \cref{eq:HAMG}.  

\subsection{Time rescaling}
\label{sec:rescal}
	Let us consider a special case with $\sigma(t)\eta(t)=1$. It follows immediately that \cref{eq:HAMG} reads as
	\begin{equation}
		\label{eq:HAMG-sigma-eta}
		\left\{
		\begin{aligned}
			{}&		X'(t)
			=
			V(t)-X(t)
			-
			\frac{1}{\eta(t)}
			\proj_{C(X(t))}
			\bigl(\eta(t)(V(t)-X(t))\bigr),\\
			{}&	\eta(t)	V'(t)
			=
			-\mu(V(t)-X(t))
			-
			\proj_{C(X(t))}
			\bigl(\eta(t)(V(t)-X(t))\bigr).
		\end{aligned}
		\right.
	\end{equation}
	Also, the parameter equation \cref{eq:HAMG-eta} becomes $		\eta'(t)
	=
	\mu-\eta(t)$ with $	\eta(0)=\eta_0>0$.
There is an explicit solution $\eta(t) = \mu+(\eta_0-\mu)e^{-t}$.
	Correspondingly, the second-order model \cref{eq:HAMG-X-2nd} becomes 
	\begin{equation}\label{eq:hamg-proj_0}
		\eta(t)X''(t) + (\mu+\eta(t))X'(t)
		+ 	\frac{\dd}{\dd t}
		\proj_{C(X(t))}(U(t))
		+ 2	\proj_{C(X(t))}(U(t)) = 0,
	\end{equation}
	where $ A(t)=e^{2t}\eta(t)/\eta_0$ and $U(t) = \eta_0(V_0-X_0)e^{-2t}$. We claim that $U(t)\to0$ as $t\to\infty$, which can be viewed as a vanishing anchor. Therefore, if we look at \eqref{eq:HAMG-X-U-X}, then we obtain \cref{eq:hamg-mog-ava-intro}, 
	and the setting $X_0=V_0$ yields a time-rescaled variant of the original \cref{eq:MOG} system.
	
	Following the argument in \cite[Section 3.2]{Luo2022d}, introducing  the  time rescaling 
	\[
	x(s)=X(t(s)),\quad v(s)=V(t(s)),\quad \hat\eta(s)=\eta(t(s)),
	\]
	with
	\begin{equation}\label{eq:res-2-din}
		\small
		t(s) = \int_{0}^{s}\lambda(r)\dd r,\quad \lambda(r) := \left\{
		\begin{aligned}
			{}&\frac{2\sqrt{\eta_{0}}}{\sqrt{\eta_{0}}r+2},&&\mathrm{~if~}\mu=0,\\
			{}&	\sqrt{\mu}\cdot\frac{(\sqrt{\eta_0}+\sqrt{\mu})e^{\sqrt{\mu}r}+\sqrt{\eta_0}-\sqrt{\mu}}{(\sqrt{\eta_0}+\sqrt{\mu})e^{\sqrt{\mu}r}+\sqrt{\mu}-\sqrt{\eta_0}},&&\mathrm{~if~}\mu>0.\end{aligned}
		\right.
	\end{equation}
Then we have $\eta(t(s))=\lambda^2(s) $ and the first-order system \cref{eq:HAMG-sigma-eta} becomes
	\begin{equation}
		\label{eq:hamg-x-v-}
		\left\{
		\begin{aligned}					\lambda(s)x'(s)={}&\lambda^2(s)\big(v(s)-x(s)\big)-\mathbf{proj}_{C(x(s))}(u(s)),\\ \lambda(s)v'(s)={}&-\mu\bigl(v(s)-x(s)\bigr)-\mathbf{proj}_{C(x(s))}(u(s)),
		\end{aligned}
		\right.
	\end{equation}
	where $u(s):=\eta_0(V_0-X_0)e^{-2 \int_{0}^{s}\lambda(r)\dd r}$. Additionally, the second-order model \cref{eq:hamg-proj_0} is transformed into 
	\[\small
	\begin{aligned}
		x''(s) + \frac{\mu+\lambda^2(s)-\lambda'(s)}{\lambda(s)}x'(s)+\frac{1}{\lambda(s)}\frac{\dd}{\dd s}
		\proj_{C(x(s))}(u(s)) +2\proj_{C(x(s))}(u(s)) = 0.
	\end{aligned}
	\]
	Therefore, if  $\mu>0$ and $\eta_0 = \mu$, then $\lambda(r) = \sqrt{\mu} $ and we obtain \cref{eq:hamh-din-intro}.
If $\mu=0$ and $\eta_0=4$, then $\lambda(r) = 2/(1+r)$ and we get \cref{eq:hamh-din-avd-intro}.


Let us look at one more time rescaling example. 	At first glance, \cref{eq:HAMG} has a first-order structure similar
to \cref{eq:HNAG-first-order}. Here, we show that for the single-objective case $m=1$, \cref{eq:HAMG} is equivalent to \cref{eq:HNAG-first-order} under proper time rescaling. Indeed, in this case \cref{eq:HAMG}  becomes
\begin{subequations}\label{eq:HAMG-1}
	\begin{empheq}[left=\empheqlbrace]{align}
		X'(t)
		={}&
		\eta(t)\sigma(t)(V(t)-X(t))
		-
		\eta(t)\sigma^{2}(t)\nabla F(X(t)),
		\label{eq:HAMG-1-X}\\
		V'(t)
		={}		&
		-\mu\sigma(t)(V(t)-X(t))
		-
		\sigma(t)\nabla F(X(t)).
		\label{eq:HAMG-1-V}
	\end{empheq}
\end{subequations}
Introduce the time transformation $s=s(t)$ by
\begin{equation}\label{eq:inverse-time-scaling}
	\frac{\dd s}{\dd t}=\sigma(t)\eta(t),\qquad s(0)=0.
\end{equation}
Since $\sigma\eta>0$, this map is strictly increasing. Write its inverse as $t=\tau(s)$ on its range, and define
\[
x(s):=X(\tau(s)),\quad v(s):=V(\tau(s)),\quad
\gamma(s):=\eta(\tau(s)),\quad \beta(s):=\sigma(\tau(s)).
\]
By the chain rule and \eqref{eq:inverse-time-scaling}, we obtain
\[
\begin{aligned}
	x'(s) =  \frac{X'(t)}{\sigma(t)\eta(t)} = V(t)-X(t)-\sigma(t)\nabla F(X(t)) = v(s)-x(s)-\beta(s)\nabla F(x(s)),
\end{aligned}
\]
and similarly,
\[
\begin{aligned}
	v'(s) = {}&V'(t)\frac{\dd t}{\dd s} = \frac{V'(t)}{\sigma(t)\eta(t)} = \frac{\mu}{\eta(t)}(X(t)-V(t))-\frac{1}{\eta(t)}\nabla F(X(t)) \\
	= {}&\frac{\mu}{\gamma(s)}(x(s)-v(s))-\frac{1}{\gamma(s)}\nabla F(x(s)).
\end{aligned}
\]
The parameter equation \cref{eq:HAMG-eta} becomes
\[
\begin{aligned}
	\gamma'(s) = \eta'(t)\frac{\dd t}{\dd s} = \frac{\eta'(t)}{\sigma(t)\eta(t)} =\mu-\eta(t) = \mu-\gamma(s).
\end{aligned}
\]
Henceforth, under the inverse time scaling \eqref{eq:inverse-time-scaling}, we recover \cref{eq:HNAG-first-order} exactly from \cref{eq:HAMG-1}.
\subsection{Energy estimate}
To derive a priori bounds on the trajectory and for further analysis, in what follows, we make the following assumption on $\sigma$. 
\begin{assum}
	\label{assum:eta-sigma} 
	The parameter $\sigma:[0,\infty)\to\R_+$ satisfies:
	\begin{itemize}
		\item[ (i)] $\sigma \in C^1$ is positive on $[0,\infty)$ and $\lim_{t\to\infty}	\int_0^t \sigma(\tau)\dd\tau=\infty$;
		\item[(ii)]  for all $t\geq0,				\theta(t):=
		\sigma'(t)
		+
		\sigma(t)
		+
		2\mu\sigma^2(t)
		-
		\frac{2+\mu\sigma(t)}{\eta(t)}\geq0$;
		\item[(iii)]   either $\sigma(t)\eta(t)= 1$ for all $t\geq0$ or $\sigma(t)\eta(t)>1$ for all $t\geq0$.
	\end{itemize}
\end{assum}
Throughout the remaining theoretical analysis, \cref{assum:eta-sigma} is imposed in addition to \cref{assum:f_j,assum:L}.
In \cref{assum:eta-sigma}, the first item is responsible for convergence rate, and the last two are used in the energy estimate. Note that there is a large class of candidates that meet \cref{assum:eta-sigma}; see \cref{tab:sigma_cond}. 
\begin{table}[H]
	\centering
	\small
	\renewcommand{\arraystretch}{1.2}
	\setlength{\tabcolsep}{3pt}
	\begin{tabular}{c c c}
		\toprule
		\multirow{2}{*}{\(\sigma(t)\)}
		& \multicolumn{2}{c}{\cref{assum:eta-sigma}} \\
		\cmidrule(lr){2-3}
		& \(\mu=0\) & \(\mu>0\) \\
		\midrule
		
		\(\dfrac{\xi(t)}{\eta(t)}\)
		& \ding{51} \((\xi(t)>1,\,\xi'\geq 2-\xi-\xi^2)\)
		& \ding{51}\ \((\xi(t)>1,\,\xi'\geq 2-\xi-\xi^2)\) \\[2ex]
		
		\(\dfrac{1}{\eta(t)}\)
		& \ding{51}
		& \ding{51} \\[1.8ex]
		
		\(c\)
		& \ding{55}
		& \ding{51}\ (\(c>1/\mu,\,\eta_0>\mu\)) \\
		
		\(c(1+\ln(1+t))\)
		& \ding{55}
		& \ding{51}\ (\(c>1/\mu,\,\eta_0>\mu\)) \\
		
		\(c(1+t)\)
		& \ding{55}
		& \ding{51}\ (\(c>1/\mu,\,\eta_0>\mu\)) \\
		
		\(ce^{t}\)
		& \ding{51}\ (\(c\geq1/\eta_0\))
		& \ding{51}\ (\(c>1/\mu,\,\eta_0>\mu\)) \\
		
		\(c(1+t)e^{t}\)
		& \ding{51}\ (\(c>1/\eta_0\))
		& \ding{51}\ (\(c>1/\mu,\,\eta_0>\mu\)) \\
		
		\bottomrule
	\end{tabular}
	\caption{A list of typical choices of \(\sigma(t)\)}
		\label{tab:sigma_cond}
\end{table}
We next derive the energy estimate. 
For later use, write $\mathbf1_m=(1,\ldots,1)^\top\in\R^m$ and define 
\begin{equation}\label{eq:R0}
	a_0:=F(X_0)+\eta_0/4\nm{V_0-X_0}^2\mathbf 1_m,
	\qquad R_0:=1+S(a_0)<+\infty.
\end{equation}
The finite property of $R_0$ comes from \cref{assum:L}.
According to the following lemma, every solution of  \cref{eq:HAMG} satisfies $X(t)\in\overline B_{R_0}$ throughout its interval of existence. 
\begin{lem}
	[\bf{Energy estimate}]
	\label{lem:f_i}
Let $(X,V)$ be a $C^1$ solution of \cref{eq:HAMG} on some interval $[0,T)$, with $T>0$. For every $0\leq t<T$, we have the energy estimate
	\begin{equation}\label{eq:H1-est-thm}
		f_j(X(t)) +	\frac12		
		\int_0^t
		\frac{\nm{X'(\tau)}^2}{\sigma^2(\tau)\eta(\tau)}\dd\tau
		\leq
		f_j(X_0)
		+
		\frac{\eta_0}{4} \nm{V_0-X_0}^2,\quad1\leq j\leq m,
	\end{equation}
	which implies
	$
	X(t)\in\mathcal L_F(a_0)\cap \overline B_{R_0}$ for all $0\leq t<T$.
\end{lem}
\begin{proof}
	Let $(X,V)$ be the given solution on $[0,T)$. 		By \cref{eq:HAMG-X-U-X}, we have
	%
	\[
	\begin{aligned}
		\dual{\nabla f_j(X),X'}={}&			\eta\sigma^{2}\big\langle \nabla f_j(X),		U-P(X,U)\big\rangle
		\leq {}\eta\sigma^{2}\big\langle P(X,U),		U-P(X,U)\big\rangle\\
		={}&\Big\langle U-\frac{X'}{	\eta\sigma^{2}},X'\Big\rangle= 		\dual{U,X'}-	\frac{1}{\sigma^2\eta}\nm{X'}^2,
	\end{aligned}
	\]
	which gives 
	\begin{equation}\label{eq:energy-est-1}
		f_j(X(t))+			
		\int_0^t
		\frac{\nm{X'(\tau)}^2}{\sigma^2(\tau)\eta(\tau)}\dd\tau\leq 		f_j(X_0)+\int_0^t\dual{U(\tau),X'(\tau)}\dd\tau.
	\end{equation}
	In view of \cref{eq:U-def-second-order}, we have
	\[
	\begin{aligned}
		\dual{U(\tau),X'(\tau)} = {}&\frac{1}{\sigma(\tau)A(\tau)} \dual{
			V_0-X_0+Z(\tau),X'(\tau)}\\
		={}&\frac{1}{\sigma(\tau)A(\tau)} \dual{
			V_0-X_0,X'(\tau)}+\frac{1}{\sigma(\tau)A(\tau)} \dual{
			Z(\tau),X'(\tau)}.
	\end{aligned}
	\]
	Applying the Cauchy--Schwarz inequality implies
	\[
	\frac{1}{\sigma(\tau)A(\tau)} \dual{
		V_0-X_0,X'(\tau)}\leq 	
	\frac{\|X'(\tau)\|^2}{2\sigma^2(\tau)\eta(\tau)}
	+
	\frac{\eta(\tau)}{2A^2(\tau)}\|V_0-X_0\|^2.
	\]
	Plugging this into \cref{eq:energy-est-1} and using \cref{lem:est-eta/A} yield that
	\begin{equation}\label{eq:est-f-key}
		f_j(X(t))+	\frac12		
		\int_0^t
		\frac{\nm{X'(\tau)}^2}{\sigma^2(\tau)\eta(\tau)}\dd\tau\leq 		f_j(X_0)+\frac{\eta_0}4\|V_0-X_0\|^2	+
		\int_0^t
		\frac{
			\dual{Z(\tau),X'(\tau)}}{\sigma(\tau)A(\tau)}
		\dd\tau .
	\end{equation}
	
	According to \cref{assum:eta-sigma} (iii), 	it holds that either $\sigma(t)\eta(t)=1$ or $\sigma(t)\eta(t)>1$. In what follows, let us prove \cref{eq:H1-est-thm} case by case. 
	First, if $\sigma(t)\eta(t)=1$ for all $t\geq0$, then $Z(t)=0$ and it is easy to obtain \cref{eq:H1-est-thm} from \cref{eq:est-f-key}.
	If $\sigma(t)\eta(t)>1$ for all $t\geq0$, then by \cref{eq:Z-A}, it is clear that $Z(0)=0$ and $Z'(t)=\frac{A(t)}{\sigma(t)\eta(t)}
	\left(1-\sigma(t)\eta(t)\right)	 X'(t)$, which gives 
	\[\small
	\begin{aligned}
		\int_0^t
		\frac{	\dual{Z(\tau),X'(\tau)}}{\sigma(\tau)A(\tau)}\dd\tau
		= 	{}&	\int_0^t
		\frac{ \eta(\tau)	\dual{Z(\tau),Z'(\tau)} }{ A^2(\tau)(1-\sigma(\tau)\eta(\tau))}
		\dd\tau 
		= {} 	\int_0^t
		\frac{ \eta(\tau)}{2 A^2(\tau)(1-\sigma(\tau)\eta(\tau))}
		\dd\nm{Z(\tau)}^2 \\
		={}&-
		\frac{ \eta(t) \nm{Z(t)}^2 }{2 A^2(t)(\sigma(t)\eta(t)-1)}
		-
		\frac12
		\int_0^t
		\frac{\eta^2(\tau)\theta(\tau)
			\nm{Z(\tau)}^2}
		{A^2(\tau)(\sigma(\tau)\eta(\tau)-1)^2}
		\dd\tau\leq0.
	\end{aligned}
	\]
	Combining this with \cref{eq:est-f-key} leads to \cref{eq:H1-est-thm}, which yields that $X(t)\in\mathcal L_{F}(a_0)$ for all $0\leq t<T$. Thanks to \cref{assum:L}, we have  $X(t)\in\overline B_{R_0}$ for all $0\leq t<T$ and  thus complete the proof. 
\end{proof}

\section{Existence of Solutions}
\label{sec:exist}
For this section, assume additionally that $\dim\mathcal H<\infty$.
More precisely, this dimensional restriction applies only to
\cref{thm:HAMG-local existence,thm:Global Existence,coro:Global Existence-res}.
The remaining sections allow infinite-dimensional Hilbert spaces
and concern any given global solution.
\subsection{Local existence}
\label{subsec:finite-local}

With $U=(V-X)/\sigma$, a solution $(X,V)$ of \cref{eq:HAMG} corresponds to a solution $(X,U)$ of \cref{eq:HAMG-X-U}, and conversely. Define a time-dependent vector field \(\mathcal G:\R_+\times\mathcal H\times\mathcal H\to \mathcal H\times\mathcal H\) by
\begin{equation}\label{eq:F}
	\mathcal G(t,\bm Z)
	:=
	\begin{bmatrix}
		\eta(t)\sigma^{2}(t)\bigl(U-P(X,U)\bigr)\\
		-\left(\sigma(t)(\mu+\eta(t))+\frac{\sigma'(t)}{\sigma(t)}\right)U
		+
		(\eta(t)\sigma(t)-1)P(X,U)
	\end{bmatrix},
\end{equation}
for all $t\geq0$ and $
\bm Z=(X,U)\in\mathcal H\times\mathcal H
$, with $P(X,U):=\proj_{C(X)}(U)$. Then \cref{eq:HAMG-X-U} is equivalent to the following non-autonomous ordinary differential equation:
\begin{equation*}
	\bm Z'(t)=\mathcal G(t,\bm Z(t)),
	\qquad
	\bm Z(0)=\bm Z_0:=(X_0,U_0)\in\mathcal H\times\mathcal H.
\end{equation*}

Below, we prove the local regularity of \(\mathcal G\), which is crucial for the subsequent local existence result.

\begin{lem}\label{lem:HAMG-continuous}
The mapping $\mathcal G$ in \cref{eq:F} is continuous on $\R_+\times\mathcal H\times\mathcal H$. Moreover, for each fixed $t\geq 0$, the mapping $\mathcal G(t,\cdot)$ is locally $1/2$-H\"older continuous with respect to $X$ and locally Lipschitz continuous with respect to $U$.
\end{lem}

\begin{proof}
Let $I\subset\R_+$ be any bounded closed interval and $\mathcal K\subset\mathcal H\times\mathcal H$ a nonempty compact (therefore bounded) subset. Set  $d_\mathcal K:=\max\{1,\sup_{\bm Z\in\mathcal K}\nm{\bm Z}\}<\infty$. 
	Since $\sigma$ is continuously differentiable and positive, there exist constants $	0<\underline\sigma\leq   \overline\sigma< \infty$ and $\widehat{\sigma}\geq0$, such that
	$
	\underline\sigma\leq \sigma(t)\leq \overline\sigma 
	$ and $\snm{\sigma'(t)}\leq \widehat{\sigma}$ for all $t\in I$. In addition, in view of \cref{eq:eta}, $\eta(t)\leq \eta_{\max}$ for all $t\geq0$. Then for any $\bm Z_1=(X_1,U_1)\in\mathcal K$ and $\bm Z_2=(X_2,U_2)\in\mathcal K$, we have
	\[\small
	\begin{aligned}
		\nm{\mathcal G(t,\bm Z_1)-\mathcal G(t,\bm Z_2)}\leq {}&	\eta(t)\sigma^{2}(t)\nm{ U_1-U_2+P(X_2,U_2)-P(X_1,U_1) }\\
		{}&\quad +\snm{ \eta(t)\sigma(t)-1}\nm{ 
P(X_1,U_1)-P(X_2,U_2)}\\
		{}&\qquad+\snm{\sigma(t)(\mu+\eta(t))+\frac{\sigma'(t)}{\sigma(t)}}\nm{		U_1-U_2}\\
		\leq{}&\left(\eta(t)\sigma^2(t)+\sigma(t)(\mu+\eta(t))+\frac{\snm{\sigma'(t)}}{\sigma(t)}\right)\nm{		U_1-U_2}\\
		{}&\quad+\left(\eta(t)\sigma^2(t)+  \eta(t)\sigma(t)+1  \right)\nm{ 
P(X_1,U_1)-		P(X_2,U_2)}\\
		\leq{}&A_1\nm{		U_1-U_2} +A_2\nm{ 
		P(X_1,U_1)-		P(X_2,U_2)},
	\end{aligned}
	\]
	with 
	\[
	A_1:= \eta_{\max}\overline{\sigma}^2+\overline{\sigma}(\mu+\eta_{\max})+\frac{\widehat{\sigma}}{\underline{\sigma}},\quad A_2: =  \eta_{\max}\overline{\sigma}^2+  \eta_{\max}\overline{\sigma}+1  .
	\]
	By \cref{eq:Proj_Lipschitz} and \cref{lem:proj_U}, it holds that
	\begin{equation*}
		\begin{aligned}
			&\|P(X_1,U_1)-		P(X_2,U_2)\|  \\
			\leq{}&
			\|P(X_1,U_1)-		P(X_1,U_2)\| 
			+\|P(X_1,U_2)-		P(X_2,U_2)\|  \\
			\leq{}&
			\|U_1-U_2\|
			+
			\sqrt{L_{d_{\mathcal K}}A_3}
			\|X_1-X_2\|^{1/2},
		\end{aligned}
	\end{equation*}
	where $A_3:=2\bigl(d_\mathcal K+M_{d_{\mathcal K}}\bigr)$.	
	Consequently, we obtain 
	\[
	\nm{\mathcal G(t,\bm Z_1)-\mathcal G(t,\bm Z_2)}\leq (A_1+A_2)\nm{U_1-U_2}	+
	A_2\sqrt{L_{d_{\mathcal K}}A_3}
	\|X_1-X_2\|^{1/2}.
	\]	
	Thus, $\mathcal G(t,\cdot)$ is $1/2$-H\"older continuous with respect to $X$ and Lipschitz continuous with respect to $U$ on $I\times \mathcal K$. 
	This concludes the proof of this lemma.
\end{proof}

According to Peano's theorem (see \cite[Theorem 2.1]{hartman2002ordinary} for instance), we obtain the local existence result.
\begin{thm}[{\bf Local existence}]\label{thm:HAMG-local existence}
Assume $\dim\mathcal H<\infty$.
	For any initial point $(X_0,V_0)\in\mathcal H\times\mathcal H$, \cref{eq:HAMG} admits a local $C^1$ solution $(X,V):[0,T_1)\to\mathcal H\times\mathcal H$ for some $T_1>0$. 
\end{thm}

	\subsection{Global existence}\label{subsec:finite-global}
	Based on the energy estimate in \cref{lem:f_i} and the classical continuation theorem for ordinary differential
	equations \cite[Corollary 3.2]{hartman2002ordinary}, we arrive at the following global
	existence result.
	\begin{thm}[{\bf Global existence}]\label{thm:Global Existence}
		Assume $\dim\mathcal H<\infty$.
		For any initial point $(X_0,V_0)\in\mathcal H\times\mathcal H$, \cref{eq:HAMG} has a global $C^1$ solution $(X,V):[0,\infty)\to\mathcal H\times\mathcal H$. In addition, we have 	$
		X(t)\in \mathcal L_F(a_0)\cap \overline B_{R_0}$ for all $t\geq0$.
	\end{thm}	
	\begin{proof}
		Invoking \cref{thm:HAMG-local existence}, \cref{eq:HAMG} admits a local
		$C^1$ solution. Let
		$ (X,V)
		$
		be a solution defined on a maximal  interval $[0,T_{\max})$, where
		$0<T_{\max}\leq\infty$. It suffices to prove that
		$T_{\max}=\infty$. 
		
		Assume, by
		contradiction, that \(T_{\max}<\infty\). Let us prove that both $X$ and $V$ are bounded on  $[0,T_{\max})$. 		
Applying \cref{lem:f_i} yields the a priori bound $\|X(t)\|\leq R_0$ for $0\leq t<T_{\max}$. Since $\proj_{C(X(t))}(U(t))$ is a convex combination of gradients at $X(t)$, it follows from \cref{eq:LR-MR} that
\[
\big\|\proj_{C(X(t))}(U(t))\big\|
\leq\max_{1\leq j\leq m}\|\nabla f_j(X(t))\|
\leq M_{R_0}<\infty,
\]
for all $0\leq t<T_{\max}$. It remains to bound $V$ on the finite interval $[0,T_{\max})$.
	Since \(\sigma\) is continuous on \([0,\infty)\), there exists
		\(0<\bar{\sigma}<\infty\) such that $		0<\sigma(t)\leq \bar{\sigma}$ for all $t\in[0,T_{\max})$. Then from \cref{eq:HAMG} we get 
		\[
		\begin{aligned}
			\|V'(t)\|
			&\leq
			\mu\sigma(t)\|V(t)-X(t)\|
			+
			\sigma(t)\left\|\proj_{C(X(t))}(U(t))\right\|  \\
			&\leq
			\mu \bar{\sigma}\|V(t)\|
			+
			\mu \bar{\sigma}R_0
			+
			\bar{\sigma}M_{R_0},
		\end{aligned}
		\]
		which implies immediately that
		\[
		\|V(t)\|
		\leq
		\|V_0\|
		+
		t(			\mu \bar{\sigma}R_0
		+
		\bar{\sigma}M_{R_0})
		+
		\mu \bar{\sigma}\int_0^t \|V(\tau)\|\dd\tau,\quad 0\leq t<T_{\max}.
		\]		
		If $\mu=0$, then $		\|V(t)\|
		\leq
		\|V_0\|
		+
		t	\bar{\sigma}M_{R_0}\leq	 
		\|V_0\|
		+
		T_{\max}	\bar{\sigma}M_{R_0}<\infty$ for all $0\leq t<T_{\max}$. If $\mu>0$, then by using Gronwall's inequality, we obtain
		\[
		\|V(t)\|
		\le
		\|V_0\|\,e^{\mu\bar{\sigma}t}
		+
		\bigl(\mu\bar{\sigma}R_0+\bar{\sigma}M_{R_0}\bigr)
		\frac{e^{\mu\bar{\sigma}t}-1}{\mu\bar{\sigma}} ,\quad 0\leq t<T_{\max}.
		\]
		Therefore, for both cases, \(V(t)\) is uniformly bounded on \([0,T_{\max})\). 
		
Since $\sigma$ is positive and continuous on $[0,T_{\max}]$, it is bounded away from zero there. Hence $U=(V-X)/\sigma$ is bounded, and the right-hand side of \cref{eq:HAMG} is bounded along the solution on $[0,T_{\max})$, which implies that $X'$ and $V'$ are bounded also. Applying the continuation principle in \cite[Corollary 3.2]{hartman2002ordinary} leads to an extended $C^1$ solution on $[0,\widehat{T})$ with $\widehat{T}>T_{\max}$. This contradiction proves $T_{\max}=\infty$ and concludes the proof.
	\end{proof}
	
	According to \cref{sec:rescal}, for the special case $\sigma\eta=1$, under the time rescaling \cref{eq:res-2-din}, we obtain \cref{eq:hamg-x-v-}, which leads to \cref{eq:hamh-din-intro} for $\mu>0$ and \cref{eq:hamh-din-avd-intro} for $\mu=0$, respectively. Moreover, by \cref{thm:Global Existence}, we claim that \cref{eq:hamg-x-v-} admits a global  $C^1$ solution $(x,v)$ on $[0,\infty)$.
\begin{coro}\label{coro:Global Existence-res} 
			Assume $\dim\mathcal H<\infty$.
	The rescaled model \cref{eq:hamg-x-v-} admits a global $C^1$ solution $(x,v):[0,\infty)\to\mathcal H\times\mathcal H$. 		 
\end{coro}	

\section{Fast Rate of Convergence}\label{sec:fast}
Here and in \cref{sec:x}, $\mathcal H$ is a real Hilbert space, possibly infinite-dimensional. We derive decay estimates for the merit function and the KKT residual along any global solution.
\subsection{Lyapunov analysis}

Recall that
$
f(X;z):=\min_{1\leq j\leq m}[f_j(X)-f_j(z)]
$; see \cref{eq:merit}.
Let $(X,V):[0,\infty)\to\mathcal H\times\mathcal H$ be a global $C^1$ solution to \cref{eq:HAMG} and for any \(z\in\mathcal H\), we introduce the Lyapunov function
\begin{equation}\label{eq:lyap}
	\mathcal E(t;z)
	:=
	f(X(t);z)+\frac{\eta(t)}{2}\|V(t)-z\|^2,
	\qquad t\geq0 .
\end{equation}
To establish the convergence rate of the Lyapunov function, we need to treat the nonsmooth term
\(f(X(t);z)\) along the trajectory. This is justified by the following lemma; see
\cite[Lemma 4.12]{Sonntag2024a}.
\begin{lem}[\cite{Sonntag2024a}]\label{lem:df}
	Let $(X,V):[0,\infty)\to\mathcal H\times\mathcal H$ be a global $C^1$ solution to \cref{eq:HAMG}. Then $t\mapsto f(X(t);z)$ is locally Lipschitz, hence locally absolutely continuous and differentiable for almost all $t>0$. In particular, for almost all $t>0$, there exists $1\leq j(t)\leq m$ such that
	\[
	\frac{\dd}{\dd t}f(X(t);z)=\left\langle\nabla f_{j(t)}(X(t)),X'(t)\right\rangle.
	\]
\end{lem}
%

This property yields that the Lyapunov function \cref{eq:lyap} is also locally absolutely continuous.
With this preparation, we derive the Lyapunov estimate.
\begin{thm}[{\bf Lyapunov estimate}] 
	\label{lem:lyap-hamg} Let $(X,V):[0,\infty)\to\mathcal H\times\mathcal H$ be a global $C^1$ solution to \cref{eq:HAMG}. Then for the Lyapunov function defined by \cref{eq:lyap} with any \(z\in\mathcal H\), we have for almost all $t\geq0$, 
	\begin{equation}\label{eq:lypn} 
		\frac{\dd}{\dd t}	\mathcal E(t;z)\leq -\sigma(t)\eta(t) \mathcal E(t;z)-\sigma^{2}(t)\eta(t)\big\|	\proj_{C(X(t))}(U(t))\big\|^2,
	\end{equation}
	where  $U=(V-X)/\sigma $. This also implies that
	\begin{equation}\label{eq:expoential inequality}	
		\begin{aligned}
			\Sigma(t)\mathcal E(t;z)+\int_{0}^{t}\sigma^{2}(\tau)\eta(\tau)\Sigma(\tau)\big\|\proj_{C(X(\tau))}(U(\tau))\big\|^2 \dd\tau
			\leq 
			\mathcal E(0;z),\quad \,t\geq0,
		\end{aligned}	
	\end{equation}
 where 
		\begin{equation}\label{eq:Sigma}
\Sigma(t)
:=
e^{\int_0^t\sigma(\tau)\eta(\tau)\dd\tau}
	=\left\{
	\begin{aligned}
		& 	1+
		\eta_0
		\int_0^t\sigma(\tau)\dd\tau, \quad && \mu = 0,\\
		& 	1+
		\frac{\eta_0}{\mu}
		\left(
		e^{\mu\int_0^t\sigma(\tau)\dd\tau}
		-1
		\right), \quad && \mu > 0.
	\end{aligned}
	\right. 
\end{equation}  
\end{thm}
	
\begin{proof}
		From \cref{eq:lyap,lem:df}, $\mathcal E(t;z)$ is differentiable for almost all $t>0$, and there exists $1\leq j(t)\leq m$ such that 
	\[
	\frac{\dd}{\dd t}	\mathcal E(t;z) 
	=  \frac{\eta'(t)}{2}\nm{V(t)-z}^2+\eta(t) \dual{V'(t),V(t)-z} + \dual{\nabla f_{j(t)}(X(t)),X'(t)},
	\]
	and plugging \cref{eq:HAMG,eq:HAMG-eta} into the right hand side gives 
	\[ 
	\begin{aligned}
		\frac{\dd}{\dd t}	\mathcal E(t;z)  
		={}&\frac{\sigma(t)\eta(t)(\mu - \eta(t))}{2}\nm{V(t)-z}^2\\
		{}&\quad +\sigma(t)\eta(t) \big\langle \mu (X(t) - V(t)) -  P(X(t),U(t)),V(t)-z\big\rangle\\
		& \quad +\sigma(t)\eta(t) \big\langle \nabla f_{j(t)}(X(t)), V(t) - X(t)- \sigma(t) P(X(t),U(t))\big\rangle.
	\end{aligned}
	\]
	Observing the identity
	\[
\begin{aligned}
	\dual{  X - V  ,V-z}
	= \frac{1 }{2}\big(\nm{X-z}^2-\nm{V-z}^2-\nm{X-V}^2\big),
\end{aligned}
\]
	we obtain
\begin{equation}\label{eq:dE-mid}
		\begin{aligned}
		\frac{\dd}{\dd t}	\mathcal E(t;z)  
		={}&-\frac{ \sigma(t)  \eta^2(t)}{2}\nm{V(t)-z}^2-\frac{\mu\sigma(t)\eta(t)}{2}\nm{X(t)-V(t)}^2\\
		{}&\quad +\frac{\mu\sigma(t)\eta(t)}{2}\nm{X(t)-z}^2-\sigma(t)\eta(t) \big\langle  P(X(t),U(t)),V(t)-z\big\rangle\\
		&\qquad +\sigma(t)\eta(t) \big\langle \nabla f_{j(t)}(X(t)), V(t) - X(t)- \sigma(t) P(X(t),U(t))\big\rangle.
	\end{aligned}
\end{equation}
	Since $\nabla f_j(X)\in C(X)$ for $1\leq j\leq m$, by  \cref{eq:opt-proj}, we have  
\[ \small
\begin{aligned}
	{}&\big\langle \nabla f_j(X)-	P(X,U) ,V-X-\sigma P(X,U)\big\rangle
	={}\sigma\big\langle \nabla f_j(X)-	P(X,U ) ,U-P(X,U )\big\rangle\leq0,
\end{aligned}
\]
which implies that 
\[
\begin{aligned}
	{}& \big\langle \nabla f_{j}(X), V - X- \sigma P(X,U)\big\rangle 
	\leq {} \big\langle    P(X,U ),V-X\big\rangle -\sigma\big\|	P(X,U )\big\|^2.
\end{aligned}
\] 
Inserting this into \cref{eq:dE-mid} implies
	\[
			\begin{aligned}
		\frac{\dd}{\dd t}	\mathcal E(t;z)  
		\leq{}&-\frac{ \sigma(t)  \eta^2(t)}{2}\nm{V(t)-z}^2 
		-\sigma^{2}(t)\eta(t) \bigl\|P(X(t),U(t) )\bigr\|^{2}\\
		{}&\quad +\frac{\mu\sigma(t)\eta(t)}{2}\nm{X(t)-z}^2-\sigma(t)\eta(t) \big\langle  P(X(t),U(t) ),X(t)-z\big\rangle.
	\end{aligned}
	\]
For $t>0$, there exists $\zeta(t)\in\Delta_m$ such that $	P(X(t),U(t) )
	=
	\sum_{i=1}^m\zeta_i(t)\nabla f_i(X(t))$. Thus, by the $\mu$-convexity of \(f_i\) (cf.\cref{assum:f_j}), we have 
	\[
\begin{aligned}
{}&\frac{\mu }{2}\nm{X(t)-z}^2-  \big\langle   P(X(t),U(t) ),X(t)-z\big\rangle 
\\
={}& \frac{\mu }{2}\nm{X(t)-z}^2-  	\sum_{i=1}^m\zeta_i(t)\dual{  \nabla f_i(X(t)),X(t)-z}\\
={}&	\sum_{i=1}^m\zeta_i(t)\left(\frac{\mu }{2}\nm{X(t)-z}^2-\dual{  \nabla f_i(X(t)),X(t)-z}\right)\leq\sum_{i=1}^m\zeta_i(t)\left[f_i(z)-f_i(X(t))\right]\\
={}&-\sum_{i=1}^m\zeta_i(t)\left[f_i(X(t))-f_i(z)\right]\leq-\min_{1\le i\le m}\left[f_i(X(t))-f_i(z)\right]=- f(X(t);z).
\end{aligned}
\]
Consequently, it follows that 
	\[
\begin{aligned}
	\frac{\dd}{\dd t}	\mathcal E(t;z)  
	\leq{}&-\sigma(t)\eta(t) 	\mathcal E(t;z) 
	-\sigma^{2}(t)\eta(t) \bigl\|P(X(t),U(t) )\bigr\|^{2},
\end{aligned}
\]
which proves \cref{eq:lypn}.  Since $\Sigma'=\sigma\eta\Sigma$, we have
\begin{equation}\label{eq:differential_ineq}
	\frac{\dd}{\dd t}[\Sigma(t)\mathcal E(t;z)]
	\leq-\sigma^2(t)\eta(t)\Sigma(t)
	\big\|P(X(t),U(t) )\big\|^2\quad\text{for almost all }t\geq0.
\end{equation}
Integrating \cref{eq:differential_ineq} and using $\Sigma(0)=1$ gives \cref{eq:expoential inequality} and completes the proof. 
\end{proof}

\subsection{Convergence rates of a merit function and the KKT residual}
Based on the Lyapunov estimate \cref{eq:expoential inequality}, we then establish the convergence rates of the merit function
$
u_0 
$ defined by \cref{eq:merit} and the KKT residual. 	According to \cref{eq:eta,eq:Sigma}, we find the relation $\eta(t)\Sigma(t) = \eta_0+\mu(\Sigma(t)-1)$.
Then by \cref{assum:eta-sigma}, we claim that
\begin{equation}\label{eq:inf-Sigma}
	\Sigma(t)\to\infty\quad \text{and}\quad \int_{0}^{t}\sigma^{2}(\tau)\eta(\tau)\Sigma(\tau)	  \dd\tau \to\infty,\quad \text{as}\,t\to\infty.
\end{equation} 
Recall the finite positive constant $R_0$ defined by \cref{eq:R0}, and the constant $N_{R_0}:=\max\{L_{R_0+1},M_{R_0}\}$ from \cref{eq:LR-MR,eq:u0-kkt}. The following theorem gives the decay rates in terms of $\Sigma(t)$. 
\begin{thm}[\bf Convergence rate]\label{thm:u0_s}
	Let $(X,V):[0,\infty)\to\mathcal H\times\mathcal H$ be a global $C^1$ solution to \cref{eq:HAMG}. 
	Then for all $t>0$, 
	\begin{equation}\label{u0_hamg}
		\begin{aligned}
			&\quad u_0(X(t)) + \frac1{\Sigma(t)} \int_{0}^{t}\sigma^{2}(\tau)\eta(\tau)\Sigma(\tau)\big\|\proj_{C(X(\tau))}(U(\tau))\big\|^2 \dd\tau\leq \frac{C_0}{\Sigma(t)},
		\end{aligned}
	\end{equation}
	where $U=(V-X)/\sigma$ and $\Sigma$ is defined by \cref{eq:Sigma} and $C_0:=	u_0(X_0)
	+
	\eta_0\|V_0\|^2
	+
	\eta_0R_0^2>0$. Moreover, for all $t>0$, we have 	the fast decay estimates on the KKT residual
	\begin{subnumcases}{\label{eq:kkt-est}}
		\label{eq:kkt-est-1}
 \big\|\proj_{C(X(t))}(0)\big\|^2 \leq   \frac{2C_0N_{R_0}}{\Sigma(t)} ,\\
 \label{eq:kkt-est-2}
 		\int_{0}^{t}\sigma^{2}(\tau)\eta(\tau)\Sigma(\tau)		\big\|\proj_{C(X(\tau))}(0)\big\|^2 \dd\tau\leq C_0.
	\end{subnumcases}
In particular, it holds that
	\begin{equation}
		\label{eq:kkt-rate}
		\inf_{0<\tau<t}		\big\|\proj_{C(X(\tau))}(0)\big\|^2 \leq   \frac{C_0+	\big\|\proj_{C(X_0)}(0)\big\|^2}{1+	\int_{0}^{t}\sigma^{2}(\tau)\eta(\tau)\Sigma(\tau)	  \dd\tau},\quad t>0.
	\end{equation}
\end{thm}

\begin{proof}
	Invoking \cref{lem:f_i}, we have
	$
	X(t)\in\mathcal L_F(a_0)$  for all $t\geq0$, where $a_0$ is defined by \cref{eq:R0}.
	Hence, by \cref{lem:L_F}, it holds that $	u_0(X(t))
	=
	\sup_{b\in F(\mathcal P_w\cap\mathcal L_F(a_0))}
	\inf_{z\in F^{-1}(b)}
	f(X(t);z)$.
%
According to \cref{eq:expoential inequality}, dropping the nonnegative squared-distance term in $\mathcal E(t;z)$, we obtain
	\begin{equation}\label{eq:f_with_integral}
		R(t;z):=f(X(t);z)
		+
		\Sigma^{-1}(t)
		\int_0^t
		\sigma^2(\tau)\eta(\tau)\Sigma(\tau)
		\|P(X(\tau),U(\tau) )\|^2
		\,\dd\tau
		\leq
		\frac{\mathcal E(0;z)}{\Sigma(t)}.
	\end{equation}
	For fixed $b$, $f(X_0;z)=\min_j(f_j(X_0)-b_j)$ is constant on $F^{-1}(b)$, so it can be taken outside the infimum. Now for any $b \in F(\mathcal P_w\cap\mathcal L_F(a_0))$, taking the infimum in \cref{eq:f_with_integral} over
	\(z\in F^{-1}(b)\), we get
\begin{equation}\label{eq:key-mid}
			\begin{aligned}
		\inf_{z\in F^{-1}(b)}R(t;z)  \leq
	\Sigma^{-1}(t) 
			\inf_{z\in F^{-1}(b)}
			\mathcal E(0;z).
	\end{aligned}
\end{equation}
By \cref{lem:S}, it follows that
	\[ 
\begin{aligned}
	{}&	\inf_{z\in F^{-1}(b)}
	\mathcal E(0;z)
	=
	\inf_{z\in F^{-1}(b)} \left(	f(X_0;z)
	+
	\frac{\eta_0}{2}\|V_0-z\|^2\right) \\
	\leq {} &	\inf_{z\in F^{-1}(b)} \left(	f(X_0;z)
	+
\eta_0\|V_0\|^2+\eta_0\|z\|^2\right)
	\\
\leq 	{}&
\eta_0\|V_0\|^2+	\inf_{z\in F^{-1}(b)}  	f(X_0;z)
	+
	\eta_0\inf_{z\in F^{-1}(b)}  \|z\|^2  
	\leq C_0.
\end{aligned}
	\] 
	Therefore, taking the supremum in \cref{eq:key-mid} over $b\in F(\mathcal P_w\cap\mathcal L_F(a_0))$ gives \cref{u0_hamg}. 		
	
	By \cref{lem:f_i}, we have $X(t)\in\overline B_{R_0}$ for all $t\geq0$, and using \cref{eq:u0-kkt}  gives \eqref{eq:kkt-est-1}.
	Since
	$
	\proj_{C(X(t))}(U(t))\in C(X(t))
	$
	and \(\proj_{C(X(t))}(0)\) is the minimum-norm element of
	\(C(X(t))\), it is clear that
	\begin{equation}\label{eq:kkt-ineq}
		\big\|\proj_{C(X(t))}(0)\big\|
		\leq
		\big\|\proj_{C(X(t))}(U(t))\big\|.
	\end{equation}
The weighted integral estimate \eqref{eq:kkt-est-2} follows directly from \cref{u0_hamg,eq:kkt-ineq}. Based on  \eqref{eq:kkt-est-2}, we obtain 
\[
\begin{aligned}
{}&	\left(1+	\int_{0}^{t}\sigma^{2}(\tau)\eta(\tau)\Sigma(\tau)	  \dd\tau\right)		\inf_{0<\tau<t}		\big\|\proj_{C(X(\tau))}(0)\big\|^2\\
	\leq{}& \big\|\proj_{C(X_0)}(0)\big\|^2+
	\int_{0}^{t}\sigma^{2}(\tau)\eta(\tau)\Sigma(\tau)		\big\|\proj_{C(X(\tau))}(0)\big\|^2 \dd\tau
		\leq \big\|\proj_{C(X_0)}(0)\big\|^2+C_0.
\end{aligned}
\]
This implies the final inequality \cref{eq:kkt-rate} and  concludes the proof of this theorem.
\end{proof}

In \cref{sec:rescal}, we have discussed \cref{eq:HAMG} with the special case $\sigma(t)\eta(t)=1$. This satisfies \cref{assum:eta-sigma} and yields a simpler model \cref{eq:HAMG-sigma-eta}, for which we have
\begin{equation}\label{eq:rate-sg-eta-1}
	u_0(X(t)) +   \int_{0}^{t}\eta^{-1}(\tau)e^{\tau-t}\big\|\proj_{C(X(\tau))}(0)\big\|^2 \dd\tau\leq C_0e^{-t},
\end{equation}
by using \cref{thm:u0_s} and the fact \cref{eq:kkt-ineq}. Moreover, under the transformation \cref{eq:res-2-din}, we obtain the rescaled model \cref{eq:hamg-x-v-}, which is related to \cref{eq:hamh-din-intro} for $\mu>0$ and \cref{eq:hamh-din-avd-intro} for $\mu=0$, respectively. In the following, we give the convergence rate of the rescaled model \cref{eq:hamg-x-v-}.
\begin{coro}
		Let $(x,v):[0,\infty)\to\mathcal H\times\mathcal H$ be a global $C^1$ solution to the rescaled model \cref{eq:hamg-x-v-} with respect to the transformation \cref{eq:res-2-din}. 
	Then for all $s>0$,  we have
\begin{equation}\label{eq:rate-sg-eta-1-res}
	u_0(x(s)) +\int_0^s \lambda^{-1}(r)\, e^{\int_{s}^{r}\lambda(w)\dd w}\, \bigl\|\proj_{C(x(r))}(0)\bigr\|^2\dd r \le C_0e^{-\int_{0}^{s}\lambda(r)\dd r},
\end{equation}
and 
\begin{equation}\label{eq:mdin-kkt}
		\big\|\proj_{C(x(s))}(0)\big\|^2 \leq \frac{2C_0N_{R_0}}{e^{\int_0^s\lambda(r)\dd r}},\quad 	\inf_{0<r<s}
		\big\|\proj_{C(x(r))}(0)\big\|^2
		\leq \frac{C_0+ \big\|\proj_{C(X_0)}(0)\big\|^2}{1+\int_0^s\lambda^{-1}(r)e^{\int_0^r\lambda(w)\dd w}\dd r},
\end{equation}
where $C_0> 0$ is defined in \cref{u0_hamg}. In particular, for $\mu=0$ and $\eta_0=4$, it follows
\[ 
\begin{aligned}
	u_0(x(s))+  \int_{0}^{s}\frac{(1+r)^3}{2(1+s)^{2}}\big\|\proj_{C(x(r))}(0)\big\|^2 \dd r\leq
C_0	(1+s)^{-2},
 \end{aligned}
\]
and
\[
\big\|\proj_{C(x(s))}(0)\big\|^2 \leq \frac{2C_0N_{R_0}}{(1+s)^2},\quad 	\inf_{0<r<s}
\big\|\proj_{C(x(r))}(0)\big\|^2
\leq \frac{8C_0 +8\big\|\proj_{C(X_0)}(0)\big\|^2}{(1+s)^4}.
\]
For $\mu>0$ and $\eta_0=\mu$, it holds
\[ 
\begin{aligned}
	u_0(x(s))+ \frac1{\sqrt{\mu}} \int_{0}^{s}e^{\sqrt{\mu}(r-s)}\big\|\proj_{C(x(r))}(0)\big\|^2 \dd r\leq
	C_0e^{-\sqrt{\mu}s},
\end{aligned}
\]
and 
\[
	\big\|\proj_{C(x(s))}(0)\big\|^2 \leq 2C_0N_{R_0}e^{-\sqrt{\mu}s},\quad 	\inf_{0<r<s}
	\big\|\proj_{C(x(r))}(0)\big\|^2
	\leq \frac{\mu C_0+\mu \big\|\proj_{C(X_0)}(0)\big\|^2}{e^{\sqrt{\mu}s}-1+\mu}.
\]
	\end{coro}
\begin{proof}
	With direct calculations, we claim that the scaling factor in \cref{eq:res-2-din} satisfies $\lambda(s)=\sqrt{\eta(t(s))}$. Moreover, we have $\Sigma(t(s))=e^{\int_0^s\lambda(r)\dd r}$ and
	\[
	\int_0^{t(s)}\sigma^2(\tau)\eta(\tau)\Sigma(\tau)\dd\tau
	=\int_0^s\lambda^{-1}(r)e^{\int_0^r\lambda(w)\dd w}\dd r.
	\]
	Plugging this into \cref{eq:kkt-est,eq:kkt-rate,eq:rate-sg-eta-1} gives \cref{eq:rate-sg-eta-1-res,eq:mdin-kkt}. 
For $\mu=0$ and $\eta_0=4$, we have $\lambda(s)=2/(1+s)$, 
and for $\mu>0$ and $\eta_0=\mu$, we have $\lambda(s)=\sqrt\mu$.
Substituting these expressions gives the stated estimates and completes the proof.
	\end{proof}

\begin{rem}
	For typical choices of $\sigma$ given in \cref{tab:sigma_cond}, we have explicit growth rates; see \cref{tab:beta_convergence_rate-}.   
	\begin{table}[H]
		\centering
		\renewcommand{\arraystretch}{1.5}
		\setlength{\tabcolsep}{3pt}
		\begin{tabular}{c c c c c}
			\toprule
			\multirow{2}{*}{\(\sigma(t)\)}
			& \multicolumn{2}{c}{\(\Sigma(t)\)}
			& \multicolumn{2}{c}{\(\sigma^2(t)\eta(t)\Sigma(t)\)} \\
			\cmidrule(lr){2-3} \cmidrule(lr){4-5}
			& \(\mu=0\) & \(\mu>0\) & \(\mu=0\) & \(\mu>0\) \\
			\midrule
			
			\(\dfrac{\xi(t)}{\eta(t)}\)
			& \(e^{\int_0^t\xi(\tau)\dd\tau}\)
			& \(e^{\int_0^t\xi(\tau)\dd\tau}\)
			& \(\xi^2(t)e^{2\int_0^t\xi(\tau)\dd\tau}\)
			& \(\xi^2(t) e^{\int_0^t\xi(\tau)\dd\tau}\)
			\\[2ex]
			
			\(\dfrac{1}{\eta(t)}\)
			& \(e^{t}\) & \(e^{t}\)
			& \(e^{2t}\) & \(e^{t}\)
			\\
			
			\(c\)
			& \(\boxed{t}\) & \( e^{\mu c t}\)
			& \(1\) & \(e^{\mu c t}\) \\
			
			\(c(1+\ln(1+t))\)
			& \(\boxed{t\ln t}\)
			& \(e^{\mu c(1+t)\ln(1+t)}\)
			& \(\ln^2 t\)
			& \(e^{\mu c(1+t)\ln(1+t)}\ln^2t \) \\
			
			\(c(1+t)\)
			& \(\boxed{t^{2}}\) & \(e^{\mu c(t+t^2/2)}\)
			& \(t^{2}\) & \(t^{2}e^{\mu c(t+t^2/2)}\) \\
			
			\(ce^{t}\)
			& \(e^{t}\) & \(e^{\mu ce^{t}}\)
			& \(e^{2t}\) & \(e^{2t+\mu ce^{t}}\) \\
			
			\(c(1+t)e^{t}\)
			& \(t e^{t}\) & \(e^{\mu cte^{t}}\)
			& \(t^2 e^{2t}\) & \(t^2 e^{2t+\mu cte^{t}}\) \\
			
			\bottomrule
		\end{tabular}
	\caption{Growth factors for $\Sigma(t)$ and $\sigma^2(t)\eta(t)\Sigma(t)$ (up to positive constant factors as $t\to\infty$); the $\mu=0$ entries in the three rows containing a boxed term are formal growth estimates only, since \cref{assum:eta-sigma} fails for these instances.}
			\label{tab:beta_convergence_rate-}
	\end{table}
\end{rem}

\section{Convergence Results of the Trajectory}\label{sec:x}
In this section, we establish the convergence of the trajectory of \cref{eq:HAMG}.  To do this, in addition to the standing \cref{assum:eta-sigma}, we need an extra restriction on $\sigma$ stated below. From \cref{tab:sigma_cond}, a typical instance is $\sigma=\xi/\eta$ with a constant $1\leq \xi\leq6$, and any time-dependent $\xi(t)\in(1,6)$ with $\xi'(t)\geq2-\xi(t)-\xi^2(t)$ is also admissible. 
\begin{assum}\label{assum:sigma-upper}
For every $t\geq0$, $\eta^2(t)\sigma(t)\leq 2\mu+6\eta(t)$.
\end{assum}
Throughout this section, \cref{assum:sigma-upper} is imposed together with \cref{assum:f_j,assum:L,assum:eta-sigma}. 
%
Let $(X,V):[0,\infty)\to\mathcal H\times\mathcal H$ be a given global $C^1$ solution of \cref{eq:HAMG}, then by \cref{lem:f_i}, $X(t)\in\overline B_{R_0}$ for all $t\geq0$. Thus $X$ is bounded and has a weak sequential cluster point. Applying the convexity gives
\begin{equation}\label{eq:trajectory-lower-bound}
f_j(X(t))\geq f_j(0)+\dual{\nabla f_j(0),X(t)}
\geq f_j(0)-R_0\|\nabla f_j(0)\|> -\infty,
\end{equation}
for all $1\leq j\leq m$.
Thus each objective is bounded below along the trajectory. Below, we first establish the monotonicity of a family of  energy functions and the convergence of the objective values. We then prove distance limits and uniqueness of weak cluster points using a
difference of squared distances, an argument developed in \cite{2025Point}.
 
To move on, let us start with the energy function
\begin{equation}\label{eq:Wi}
	\mathcal W_j(t): =
	 f_j(X(t))+\frac{\eta(t)}{2}\nm{V(t)-X(t)}^2,\quad 1\leq j\leq m.
\end{equation}
We shall establish the nonincreasing property of $\mathcal W_j$.
\begin{lem}\label{lem:basic_dissipation}  
	Let $(X,V):[0,\infty)\to\mathcal H\times\mathcal H$ be a global $C^1$ solution  to \cref{eq:HAMG}.
	For every $1\leq j\leq m$, the energy function \(\mathcal W_j\) defined by \cref{eq:Wi} is nonincreasing on 
	\([0,\infty)\).
\end{lem}  
\begin{proof}
	Subtracting the first equation of
	\cref{eq:HAMG} from the second one, we obtain
\begin{equation}\label{eq:V'-X'}
	(V-X)'
	={}
	-\sigma^2(\mu+\eta)U+
	\sigma(\eta\sigma-1)
	P(X,U).
\end{equation}
Then by \cref{eq:HAMG-eta} and \eqref{eq:HAMG-X-U-X}, a direct computation yields that
	\[\small
\begin{aligned}
	\frac{\dd}{\dd t}\mathcal W_j
	={}&\dual{\nabla f_j(X),X'} +\frac{\eta'}{2}\nm{V-X}^2+\eta\dual{V-X,V'-X'} \\
	={}&		\eta\sigma^2
	\big\langle
	\nabla f_j(X),
	U-P(X,U)
	\big\rangle-
	\frac{\eta\sigma^3}2
	(\mu+3\eta)
	\|U\|^2+
	\eta\sigma^2
	(\eta\sigma-1)
	\big\langle
	U,
	P(X,U)
	\big\rangle.
\end{aligned}
\]
	Invoking \cref{eq:opt-proj} gives
	\[
\begin{aligned}
	{}&	 
	\big\langle
	\nabla f_j(X),
	U-P(X,U)
	\big\rangle+
	(\eta\sigma-1)
	\big\langle
	U,
	P(X,U)
	\big\rangle \\
	\leq{} &
	\eta\sigma
	\big\langle
	U,P(X,U)
	\big\rangle
	-
	\big\|P(X,U)\big\|^2
	\leq{}\frac{	\eta^2\sigma^2 }{4}\nm{U}^2.
\end{aligned}
\]
	Therefore, we obtain
	\[
\begin{aligned}
	\frac{\dd}{\dd t}\mathcal W_j
	\leq{}&
	\frac{	\eta^3\sigma^4 }{4}\nm{U}^2-
	\frac{\eta\sigma^3}2
	(\mu+3\eta)\|U\|^2
	={}
	-\frac{	\eta\sigma^3 }{4}
	\left[		2
	\mu+6\eta-\eta^2\sigma
	\right]
	\|U\|^2\leq0,
\end{aligned}
\]
	where the last step follows from \cref{assum:sigma-upper}. 
	This completes the proof.
\end{proof}
%
%
%
\begin{lem}\label{lem:limit-fi-HAMG}
	Let $(X,V):[0,\infty)\to\mathcal H\times\mathcal H$ be a global $C^1$ solution  to \cref{eq:HAMG}. Then we have the following.
	\begin{itemize}
		\item[(i)] $D:=V-X\in C^1([0,\infty);\mathcal H)\cap L^2(0,\infty;\mathcal H)$ and $\lim_{t\to\infty}D(t)=0$;
		\item[(ii)] $\nm{V(t)}\leq R_0+\nm{V_0-X_0}+\sqrt{5C_0/\eta_0}$ for all $t\geq0$;
		\item[(iii)] for $1\leq j\leq m$, there is a finite real number  $\ell_j\in\R$ such that $		\lim_{t\to\infty}f_j(X(t))=\lim_{t\to\infty}\mathcal W_j(t)=\ell_j$;
		\item[(iv)] the set $\mathcal X:=\{z\in\mathcal H:f_j(z)\leq\ell_j,\,1\leq j\leq m\} \subset \mathcal P_w$ is nonempty and contains all weak cluster points of $X$, and for any $z\in\mathcal X$, we have $	f_j(X(t))-f_j(z)\geq -\eta(t)/2\|D(t)\|^2$ for all $t\geq0$.
	\end{itemize} 
\end{lem}
\begin{proof}
By \cref{eq:eta,eq:Sigma}, we have $\eta(t)\Sigma(t)=\eta_0e^{\mu\int_0^t\sigma(s)\dd s}\geq\eta_0$. The integral estimate in \cref{u0_hamg}, together with the fact $u_0(X(t))\geq0$, leads to
\begin{equation}\label{eq:sigmaP-L2}\small
\int_0^t\sigma^2(\tau)\|P(X(\tau),U(\tau))\|^2\dd\tau
\leq\frac1{\eta_0}\int_0^t\sigma^2(\tau)\eta(\tau)\Sigma(\tau)\|P(X(\tau),U(\tau))\|^2\dd\tau
\leq\frac{C_0}{\eta_0}.
\end{equation}
Set $a(t):=\eta(t)\sigma(t)$ and $b(t):=\sigma(t)(\mu+\eta(t))$. By \cref{assum:eta-sigma,assum:sigma-upper}, we find that
\[
1\leq a\leq6+\frac{2\mu}{\eta},\qquad
b=a\left(1+\frac{\mu}{\eta}\right)\geq1,
\qquad
\frac{(a-1)^2}{b}
\leq\frac{a-1}{1+\mu/\eta}
\leq\frac{5+2\mu/\eta}{1+\mu/\eta}\leq5.
\]
Using \cref{eq:V'-X'}, the relation $D=\sigma  U$ and Young's inequality, we obtain,  
\begin{align*}
\frac{\dd}{\dd t}\|D\|^2
&=-2b\|D\|^2+2\sigma(a-1)\langle P(X,U),D\rangle\\
&\leq-b\|D\|^2+\frac{(a-1)^2}{b}\sigma^2\|P(X,U)\|^2
\leq-\|D\|^2+5\sigma^2\|P(X,U)\|^2.
\end{align*}
The nonnegative term is continuous and integrable by \cref{eq:sigmaP-L2}. Applying \cref{lem:h-limit} gives $\lim_{t\to\infty}D(t)=0$. Integrating the same differential inequality yields
\[
\begin{aligned}
	\|D(t)\|^2+\int_0^t\|D(\tau)\|^2\dd\tau
	\leq{}&\|V_0-X_0\|^2+5\int_0^t\sigma^2(\tau)\|P(X(\tau),U(\tau))\|^2\dd\tau\\
	\leq{}&\|V_0-X_0\|^2+\frac{5C_0}{\eta_0}.
\end{aligned}
\]
Letting $t\to\infty$ proves $D\in L^2(0,\infty;\mathcal H)$. Since $V=X+D$ and $\|X(t)\|\leq R_0$, the bound in (ii) follows.

According to \cref{lem:basic_dissipation}, $\mathcal W_j$ is nonincreasing and bounded below by the finite lower bound for $f_j(X(t))$ in \cref{eq:trajectory-lower-bound}. Hence it  converges to a finite limit $\ell_j\in\R$. In particular, we have
\begin{equation}\label{eq:lb-1}
	f_j(X(t))+\frac{\eta(t)}2\|D(t)\|^2
	=\mathcal W_j(t)\geq\ell_j,\quad \forall\,t\geq0.
\end{equation}
%
%
%
Since $0<\eta(t)\leq\eta_{\max}$, we have $0\leq\mathcal W_j(t)-f_j(X(t))=\eta(t)/2\|D(t)\|^2\longrightarrow0$ as $t\to\infty$,
which gives $f_j(X(t))\to\ell_j$ as well.

For any weak cluster point $x_\infty\in\mathcal H$  of $X(t)$, there exists $\{X(t_k)\}_{k\in\mathbb N}$ with $t_k\to\infty$ as $k\to\infty$, such that $X(t_k)\rightharpoonup  x_\infty $ as $k\to\infty$. Thanks to the convexity and continuity, each $f_j$ is weakly lower semicontinuous, and therefore $f_j(x_\infty)\leq\liminf_{k\to\infty}f_j(X(t_k))=\ell_j$, which means $\mathcal X$ contains all weak cluster points of the trajectory $X$. Boundedness of $X$ and reflexivity of $\mathcal H$ ensure that at least one such weak cluster point exists, so $\mathcal X$ is nonempty. It remains to show $\mathcal X\subset\mathcal P_w$. Indeed, otherwise some $x\in\mathcal X$ would admit $y\in\mathcal H$ with $F(y)\prec F(x)$. Since $m$ is finite, we have $\delta:=\min_{1\leq j\leq m}\bigl(\ell_j-f_j(y)\bigr)>0$.
The definition of $u_0$ and convergence of the objective values imply
\[
\liminf_{t\to\infty}u_0(X(t))
\geq\lim_{t\to\infty}\min_{1\leq j\leq m}\bigl(f_j(X(t))-f_j(y)\bigr)=\delta,
\]
contradicting $u_0(X(t))\to0$ as $t\to\infty $ from \cref{eq:inf-Sigma,thm:u0_s}. Finally, the estimate $	f_j(X(t))-f_j(z)\geq -\eta(t)/2\|D(t)\|^2$ follows from \cref{eq:lb-1} immediately. This concludes the proof.
\end{proof}
%
%
%

Based on the previous results, we are ready to state the main result of this section.
\begin{theorem}[\bf Weak and strong convergence]\label{thm:sc}
Under Assumptions~\ref{assum:f_j}, \ref{assum:L}, \ref{assum:eta-sigma} and \ref{assum:sigma-upper}, for every global $C^1$ solution $(X,V):[0,\infty)\to\mathcal H\times\mathcal H$ to \cref{eq:HAMG}, both $X$ and $V$ converge weakly to a weak Pareto point $x_\infty\in\mathcal P_w$. The convergence is strong if $\mu>0$. 
\end{theorem}
\begin{proof}
Fix $z \in\mathcal X$. For each $t>0$, there is $\lambda(t)\in\Delta_m$ such that $P(X(t),U(t))=\sum_{j=1}^m\lambda_j(t)\nabla f_j(X(t))$.
Applying the $\mu$-convexity and \cref{lem:limit-fi-HAMG} (iv) gives
\[
\langle P(X,U),X-z \rangle
\geq\sum_{j=1}^m\lambda_j\bigl(f_j(X)-f_j(z )\bigr)
+\frac\mu2\|X-z \|^2
\geq -\frac{\eta}2\|V-X\|^2+\frac\mu2\|X-z \|^2.
\] 
Set $\psi_{z }(t):=\|V(t)-z \|^2/2$. From the second equation in 	\cref{eq:HAMG}, we obtain $V'=-\mu\sigma D-\sigma P(X,U)$ and 
\begin{align*}
	\psi_{z }'
	={}&-\mu\sigma\langle D,V-z \rangle-\sigma\langle P(X,U),X-z \rangle-\sigma\langle P(X,U),D\rangle\\
	={}&-\frac{\mu\sigma}2\left(\|V-z \|^2+\|D\|^2-\|X-z \|^2\right)-\sigma\langle P(X,U),X-z \rangle-\sigma\langle P(X,U),D\rangle\\
	\leq{}&-\mu\sigma\psi_{z }+\frac{\sigma(\eta-\mu)}2\|D\|^2+\sigma\|P(X,U)\|\|D\|.
\end{align*} 
The additional upper bound in \cref{assum:sigma-upper} yields $\sigma(\eta-\mu)\leq6$.
Thus, we have
\begin{equation}\label{eq:psi}
	\psi_{z }'+\mu\sigma\psi_z\leq3\|D\|^2+\sigma\|P(X,U)\|\|D\| \leq4\|D\|^2+\sigma^2\|P(X,U)\|^2:=g,
\end{equation}
By \cref{eq:sigmaP-L2} and \cref{lem:limit-fi-HAMG} (i), the right hand side $g\in C[0,\infty)\cap L^1(0,\infty)$, and invoking 
\cref{lem:h-limit} implies $\lim_{t\to\infty}\psi_z(t)$ exists. Since $X,V$ are bounded and $\lim_{t\to\infty}D(t)=0$ (cf.\cref{lem:f_i,lem:limit-fi-HAMG}), it follows that
\begin{equation}\label{eq:X-z}\small
	\left|\|X(t)-z \|^2-\|V(t)-z \|^2\right|
	\leq\|D(t)\|\bigl(\|X(t)-z \|+\|V(t)-z \|\bigr)\longrightarrow0\quad \text{as } t\to\infty.
\end{equation}
Hence $\|X(t)-z \|^2$ has also a finite limit for every $z \in\mathcal X$. Let $z_1,z_2\in\mathcal X$ be two weak cluster points of $X$. Consider the difference function
\[
h(t):=\|X(t)-z_1\|^2-\|X(t)-z_2\|^2
=2\langle X(t),z_2-z_1\rangle+\|z_1\|^2-\|z_2\|^2,
\]
which has a finite limit as $t\to\infty$. Take two sequences $t_k\to\infty$ and $s_k\to\infty$ such that  $X(t_k)\rightharpoonup z_1$ and $X(s_k)\rightharpoonup z_2$. Then it follows that 
\[
\lim_{k\to\infty}h(t_k)=-\|z_1-z_2\|^2,
\qquad \lim_{k\to\infty}h(s_k)=\|z_1-z_2\|^2.
\]
These limits coincide, proving $z_1=z_2=:x_\infty$. Since $X$ is bounded, for any sequence $t_k\to\infty$, $X(t_k)$ is bounded and has $x_\infty$ as its unique weak sequential cluster point. Hence by \cite[Lemma 2.38]{BauschkeCombettes2011} we conclude that $X(t_k)\rightharpoonup x_\infty$. This, together with the arbitrariness of the sequence $t_k$, implies that $X(t)\rightharpoonup x_\infty\in\mathcal X\subset\mathcal P_w$ as $t\to\infty$. As a result we have the weak convergence $V(t)\rightharpoonup x_\infty$ as $t\to\infty$.

For $\mu>0$, we have $\int_{0}^{\infty}\mu\sigma(t)\dd t = \infty$, and applying \cref{lem:h-limit} to \cref{eq:psi} leads to $V(t)\to x_\infty$ as $t\to\infty$. Plugging this into \cref{eq:X-z} gives the strong convergence $X(t)\to x_\infty$ as $t\to\infty$. This finishes the proof.
\end{proof}

In view of the rescaling \cref{eq:res-2-din}, it is clear that $t(s)\to\infty$ as $s\to\infty$. Hence, we obtain the weak convergence of the rescaled model \cref{eq:hamg-x-v-}, which leads to \cref{eq:hamh-din-intro} for $\mu>0$ and \cref{eq:hamh-din-avd-intro} for $\mu=0$, respectively. 
\begin{coro}
	Under Assumptions~\ref{assum:f_j}, \ref{assum:L}, \ref{assum:eta-sigma} and \ref{assum:sigma-upper}, for every global $C^1$ solution $(x,v):[0,\infty)\to\mathcal H\times\mathcal H$ to the rescaled model \cref{eq:hamg-x-v-}, both $x$ and $v $ converge weakly to a weak Pareto point. The convergence is strong if $\mu>0$.
\end{coro}
 
%
\section{Numerical Tests}\label{sec:tests}
In this section, we investigate the convergence behavior of \cref{eq:HAMG} for solving strongly convex ($\mu>0$), convex ($\mu=0$) and nonconvex multiobjective problems, with regard to the fast rate of the KKT residual, the oscillation-reduction effect and the quality of the approximated Pareto front.
%
%
%

\subsection{Convergence behavior with different parameters} 
\begin{figure}[H]
	\centering
	\begin{minipage}{0.7\textwidth}
		\centering
		\includegraphics[width=\linewidth]{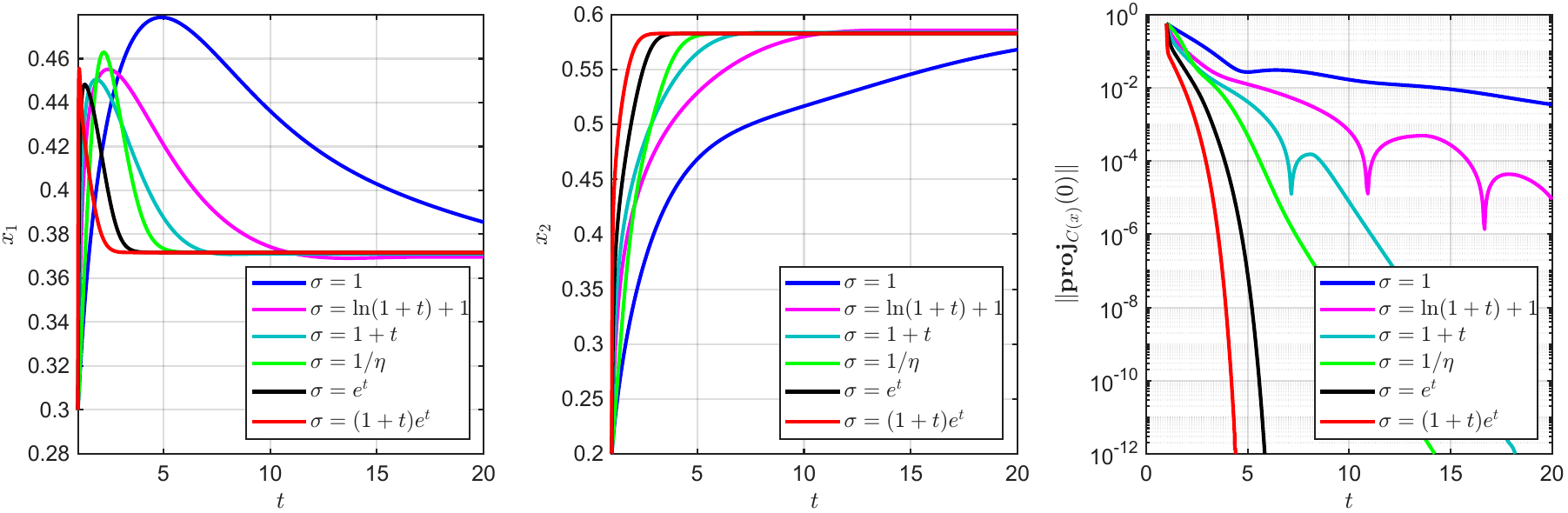}
	\end{minipage}
	\caption{Numerical result of \cref{eg:beta_sc}. Left: the first component. Middle: the second component. Right: the KKT residual.}
	\label{fig:eta-sc}
\end{figure}

\begin{figure}[H]
	\centering
	\begin{minipage}{0.7\textwidth}
		\centering
		\includegraphics[width=\linewidth]{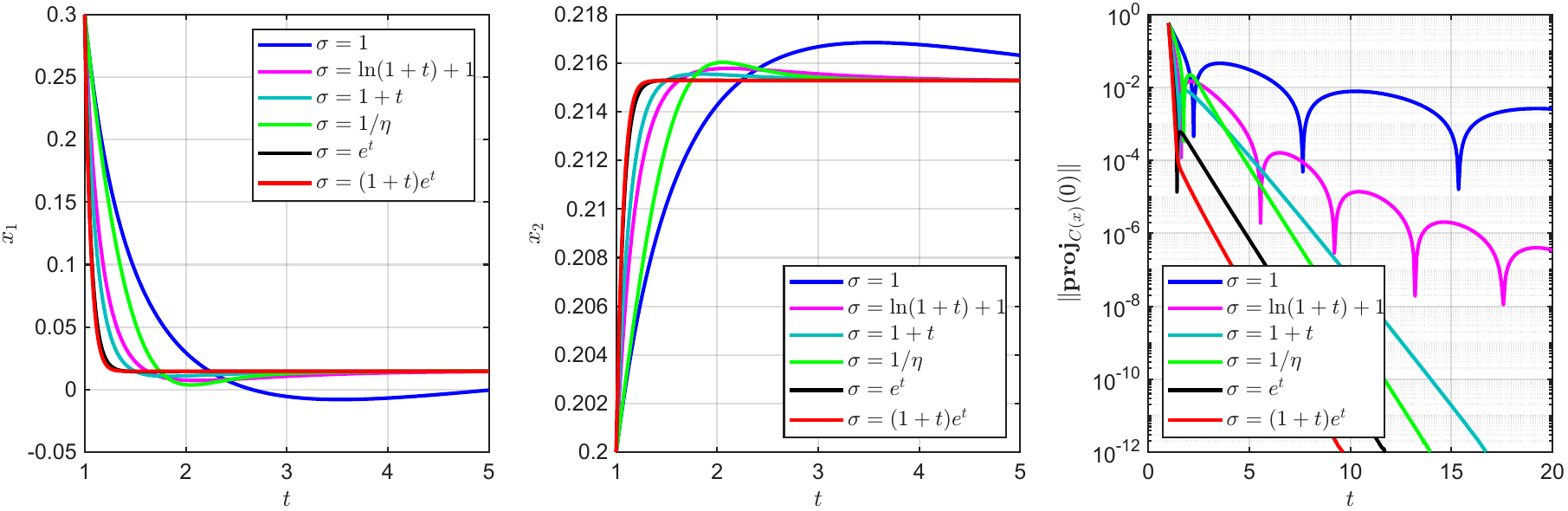}
	\end{minipage}
	\caption{Numerical result of \cref{eg:beta_c}. Left: the first component. Middle: the second component. Right: the KKT residual.}
		\label{fig:eta-c}
\end{figure}

\begin{figure}[H]
	\centering
	\begin{minipage}{0.7\textwidth}
		\centering
		\includegraphics[width=\linewidth]{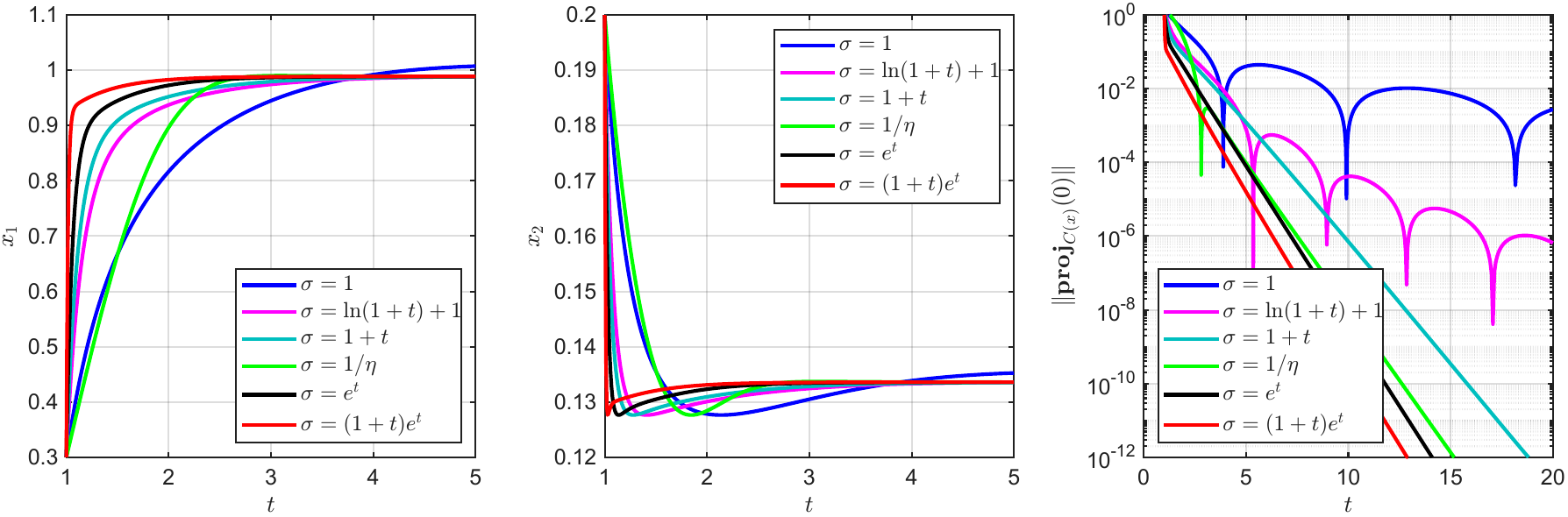}
	\end{minipage}
	\caption{Numerical result of \cref{eg:beta_nc}. Left: the first component. Middle: the second component. Right: the KKT residual.}
		\label{fig:eta-nc}
\end{figure}

This parameter study explores the behavior of HAMG with different parameters. 
In \cref{eg:beta_sc,eg:beta_c,eg:beta_nc}, we test several choices of \(\sigma(t)\) listed in \cref{tab:beta_convergence_rate-}, including choices outside \cref{assum:eta-sigma,assum:sigma-upper}. For each example, we set $\eta_0=1$ and generate the initial point $(X_0,V_0)$ randomly.  The trajectory $X(t)$ and the KKT residual $\|\proj_{C(X(t))}(0)\|$ are shown in \cref{fig:eta-sc,fig:eta-c,fig:eta-nc}. For $\sigma(t)= 1$ and $\sigma(t)=1+\ln(1+t)$, the HAMG flow converges more slowly and exhibits oscillations. For the other choices, $\sigma(t)= 1/\eta(t),\,\sigma(t)= 1+t,\,\sigma(t)= e^t$ and $\sigma(t)= (1+t)e^t$, we observe faster decay of the KKT residual. 
\begin{eg}\label{eg:beta_sc}
Consider two quadratic objective functions
	\begin{equation*}
		f_j(x) = \frac{\mu}{2} \|x\|^2 + \frac{1}{2} \|A_j x - b_j\|^2, \,x\in\R^n,\,\quad j = 1, 2,
	\end{equation*}
	where \(A_j \in \mathbb{R}^{p \times n}\) and \(b_j \in \mathbb{R}^p\).  
Here, we set \(n = 2\), \(p = 2\), and \(\mu = 0.05\) is a common strong-convexity parameter.  
	The matrices \(A_j\) and vectors \(b_j\) are generated from the uniform distribution on \([0,1]\).  
\end{eg}

\begin{eg}\label{eg:beta_c}
	Consider two convex objectives
\[
			f_1(x) = (x_1 - 2)^4 + (x_2 - 2)^8, \quad
f_2(x) = x_1^2 + x_2^2,\quad x = (x_1,x_2)\in\R^2.
\]
\end{eg}

\begin{eg}\label{eg:beta_nc}
	Consider the following nonconvex problem
\[
			f_1(x) = x_1^4 + x_2^4 - x_1^2 + x_2^2 - 10x_1x_2 + 0.25x_1 + 20, \quad
f_2(x) = (x_1 - 1)^2 + x_2^2,\quad x = (x_1,x_2)\in\R^2.
\]
\end{eg}
%


\subsection{Comparison with existing dynamical systems}
  
Let us compare \cref{eq:HAMG} with three baseline dynamical systems: \cref{eq:MOG} \cite{attouch_multiibjective_2015}, \cref{eq:AMG} \cite{Luo2025a}, and \cref{eq:MAVD} \cite{Sonntag2024a}. 
For HAMG, we take $\eta_0=1$ and set $\sigma(t)=e^t$ for $\mu=0$ and $\sigma(t)=6/\eta(t)$ for $\mu>0$.  
For AMG, we set $\gamma_0 = 1$, and for MAVD we choose $\alpha=3$. If the problem is strongly convex, then we use a common strong-convexity parameter $\mu$ for HAMG and AMG, and  for convex and nonconvex problems, we take $\mu=0$ for  HAMG and AMG. For dynamical models, we use the same initial points sampled uniformly from $[-R,R]^n$ with $R\in\{2,50\}$, and the number of initial points is \(N=200\). The systems are integrated using MATLAB's built-in solver \texttt{ode15s}, with a common integration
interval in dynamical time $[t_0,T_{\max}]$ where $0<t_0<T_{\max}<\infty$.  

To measure the quality of the approximate weak Pareto front, we introduce some quantities. Denote by $\mathcal A=\{\rm HAMG, MOG,AMG,MAVD\}$ the set of dynamical systems. For each $a\in\mathcal A$, given the initial points $\{X_0^k\in\R^n:\,1\leq k\leq N\}$, the output trajectories are $\{X_a^1(t),\cdots,X_a^N(t)\}$ and the corresponding end points are $\{X^k_a(T_{\max}):\,1\leq k\leq N\}$. Let us introduce $\Omega:=\cup_{a\in\mathcal A}\Omega_a$, where
\[
\begin{aligned}
\Omega_a:=\bigl\{F(X^k_a(T_{\max}))\in\R^m:\;&1\leq k\leq N,\ \text{there is no }1\leq\ell\leq N\\
&\text{such that }F(X^\ell_a(T_{\max}))\prec F(X^k_a(T_{\max}))\bigr\}.
\end{aligned}
\]
Here $\prec$ denotes strict inequality in every component, so $\Omega_a$ is the weakly nondominated subset of the sampled objective vectors and serves as an approximate weak Pareto front. Denote by ${\bf NP}(a)=\snm{\Omega_a}$ the number of distinct weakly nondominated objective vectors. Consistently, define the empirical reference weak Pareto front as $\Omega_{\rm ref}:=\{Y\in\Omega:\,\text{there is no }Z\in\Omega\text{ such that }Z\prec Y\}$. Once we get $\Omega_{\rm ref}$, define the purity of each $a\in\mathcal A$ by $	\operatorname{\bf Purity}(a):=
	|\Omega_a\cap\Omega_{\mathrm{ref}}|/
	|\Omega_a|$. A larger purity value indicates that a larger proportion of the distinct objective vectors in $\Omega_a$ survives weak-nondominance filtering of the pooled outputs.
	To further evaluate the quality of $\Omega_a$
	we use the hypervolume (HV) \cite{zitzler2000comparison}
	and the modified inverted generational distance (IGD$^+$)
	\cite{ishibuchi2015modified}:
	\[
	{\bf HV}(a): = {\rm volume}\left(		\cup_{y\in\Omega_a}[y,p]\right),\quad {\bf IGD^+}(a): = 	\frac{1}{|\Omega_{\mathrm{ref}}|}
	\sum_{z\in\Omega_{\mathrm{ref}}}
	\min_{y\in\Omega_a}
	d^{+}(y,z),
	\]
	where \(p\in\mathbb{R}^{m}\) is a common reference point satisfying
\(y\preceq p\) for all
\(y\in\Omega\), and $
d^{+}(y,z)
=\nm{(y-z)_+} 
$
denotes the Euclidean norm of the nonnegative part of $y-z$. 	A larger HV value yields better dominance and coverage while a smaller IGD\(^{+}\) value implies a better approximation to $\Omega_{\mathrm{ref}}$. 
Moreover, we introduce two more quantities to measure the speed of convergence. Specifically, for each $a\in\mathcal A$, with the initial points $\{X_0^k\in\R^n:\,1\leq k\leq N\}$, define two averaged residuals:
\begin{equation}\label{eq:kkt}
	{\rm Gap}_a(t): = \frac{1}{N}\sum_{k=1}^{N}\operatorname{dist}
	\bigl(
	F(X_a^{k}(t)),
	\Omega_{\mathrm{ref}}
	\bigr),\quad {\rm KKT}_a(t): = \frac{1}{N}\sum_{k=1}^{N}	\big\|
	\proj_{C(X_a^{k}(t))}(0)
	\big\|.
\end{equation}
Based on this, we define the smallest settling times after which the averaged objective-space distance and KKT residual remain at most $\varepsilon$:
\[
\begin{aligned}
	T_{\rm PF}(a): = {}& \inf
	\left\{
	t\in[t_0,T_{\max}] :
{\rm Gap}_a(s)\leq\varepsilon
	\ \text{ for all }t\leq s\leq T_{\max}
	\right\},\\
	T_{\rm KKT}(a): ={}&  \inf
	\big\{
	t\in[t_0,T_{\max}] :
{\rm KKT}_a(s)
	\leq
	\varepsilon
	\ \text{ for all }t\leq s\leq T_{\max}
	\big\},
\end{aligned}
\]
where $\varepsilon=10^{-6}$ is the given tolerance. 
If the corresponding residual does not
satisfy the settling criterion $\varepsilon$ by $T_{\max}$, then we set $\inf\emptyset:=\infty$.
\subsubsection{An imbalanced example}
Let us test the performance on an imbalanced problem, which has ill-conditioned components or components on substantially different scales.
\begin{eg}\label{eg:3}
Consider two quadratic objective functions
\[
		f_j(x)=\frac{1}{2}(x-b^j)^{\top}Q_j(x-b^j),\quad 	Q_1=\begin{bmatrix}
			2 & 0 \\ 
			0 & 1/2
		\end{bmatrix},\quad 
		Q_2=\begin{bmatrix}
			\epsilon_1 & 0 \\ 
			0 & \epsilon_2 
		\end{bmatrix},
\]
where $b^j\in\R^2$ and $j = 1,2$. 
	We are interested in two cases. One is that $\epsilon_1\ll\epsilon_2$, which means $f_2$ itself is ill-conditioned. The second is that $\epsilon_1=\epsilon_2\ll 1$, then both $f_1$ and $f_2$ are well-conditioned but the global condition number is $2/\epsilon_1$.
\end{eg} 
\subsubsection*{Case 1: $\epsilon_1=10^{-5},\,\epsilon_2=2$}
\begin{figure}[H]
	\centering
	\begin{minipage}{0.5\textwidth}
		\centering
		\includegraphics[width=\linewidth]{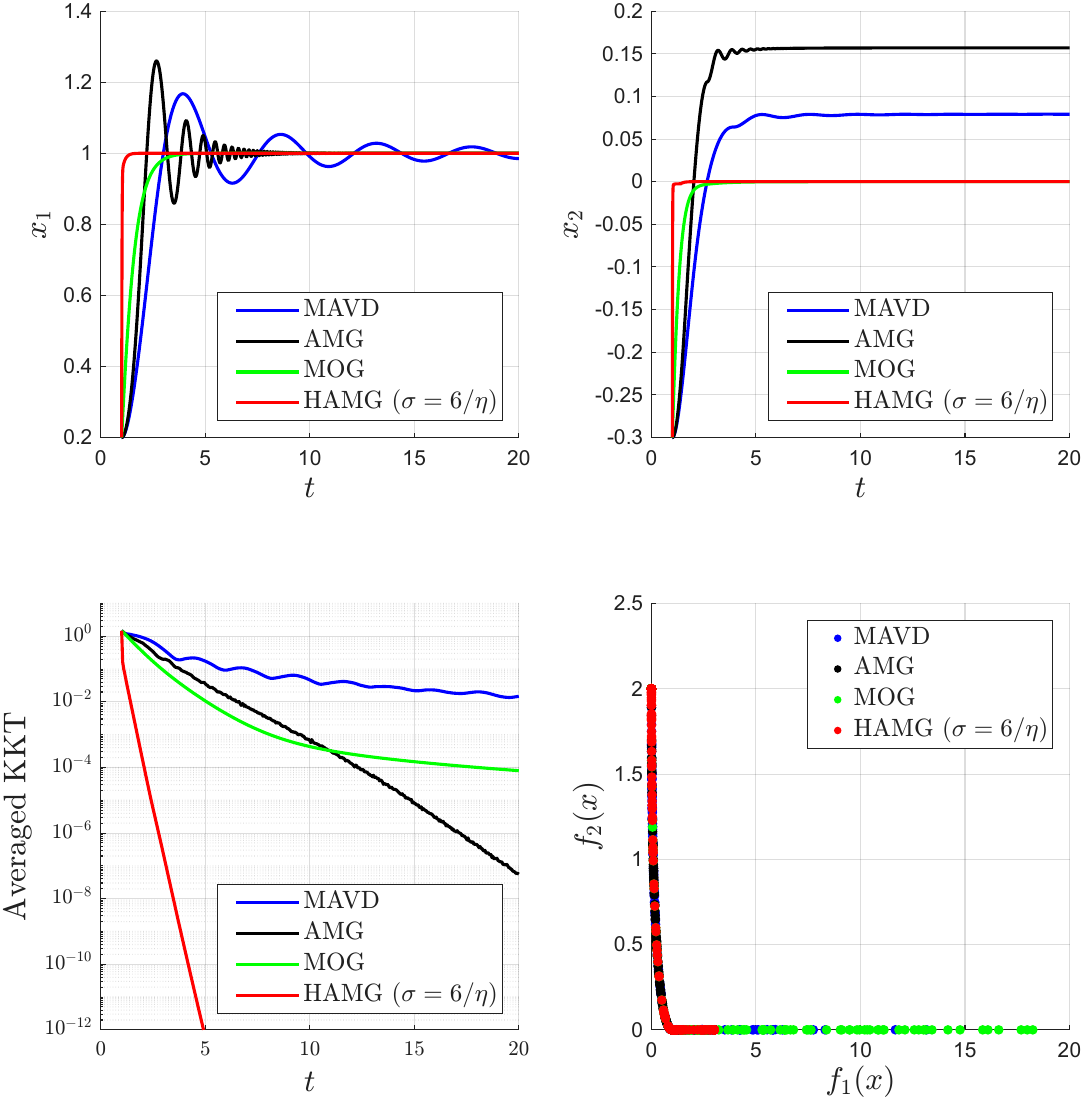}
	\end{minipage}
	\caption{Numerical results of \cref{eg:3} with $\epsilon_1=10^{-5}$ and $\epsilon_2=2$. The initial range is $[-2,2]^2$. Top line: the components of a sample trajectory. Bottom line: the averaged KKT residual (left) and  the approximate Pareto front (right).}
	\label{fig:img4}
\end{figure}

In this case, a common strong convexity parameter is $\mu=10^{-5}$. The numerical results are summarized in \cref{tab:Partially imbalanced example} and \cref{fig:img4}. As we can see, in both initial ranges, HAMG has the shortest settling times to ensure the residuals enter a local region with $\varepsilon=10^{-6}$. AMG
yields the best HV and IGD$^+$ values, while HAMG achieves high purity and competitive front
approximations.  Also, it is observed from \cref{fig:img4} that the trajectory and the averaged KKT residual (cf.\cref{eq:kkt}) of HAMG exhibit better decay rates, compared to the heavy oscillation or slow convergence behaviors of AMG and MAVD.
\begin{table}[H]
	\centering
	\caption{Performances of different dynamical models for \cref{eg:3} with $\epsilon_1=10^{-5}$ and $\epsilon_2=2$. ``$\infty$'' indicates that the corresponding residual does not
		satisfy the settling criterion by $T_{\max} = 20$.}
	\label{tab:Partially imbalanced example}
	\small
	\setlength{\tabcolsep}{10pt}
	\renewcommand{\arraystretch}{1.08}
	\resizebox{\textwidth}{!}{
		\begin{tabular}{lcccccc}
			\toprule
			Model
			& $T_{\mathrm{PF}}$
			& $T_{\mathrm{KKT}}$
			& NP
			& Purity $\uparrow$
			& HV $\uparrow$
			& IGD$^{+}$ $\downarrow$ \\
			\midrule
			
			\multicolumn{7}{c}{Initial range $[-2,2]^2$} \\
			\midrule
			
			MAVD
			& $\infty$
			& $\infty$
			& 125
			& 0.76
			& 11.65
			& $2.99\times10^{-3}$ \\
			
			AMG
			& $\infty$
			& 18.02
			& 193
			& 0.99
			& \textbf{11.66}
			& $\mathbf{1.07\times10^{-3}}$ \\
			
			MOG
			& $\infty$
			& $\infty$
			& 115
			& 0.97
			& 11.64
			& $5.37\times10^{-3}$ \\
			
			HAMG ($\sigma=6/\eta$)
			& \textbf{19.97}
			& \textbf{2.77}
			& \textbf{199}
			& \textbf{1.00}
			& 11.64
			& $3.91\times10^{-3}$ \\
			
			\midrule
			\multicolumn{7}{c}{Initial range $[-50,50]^2$} \\
			\midrule
			
			MAVD
			& $\infty$
			& $\infty$
			& 34
			& 0.32
			& 24.32
			& $9.79\times10^{-3}$ \\
			
			AMG
			& $\infty$
			& $\infty$
			& 141
			& 0.99
			& \textbf{24.40}
			& $\mathbf{1.20\times10^{-3}}$ \\
			
			MOG
			& $\infty$
			& $\infty$
			& 15
			& 0.60
			& 22.76
			& $1.70\times10^{-1}$ \\
			
			HAMG ($\sigma=6/\eta$)
			& \textbf{20.00}
			& \textbf{3.35}
			& \textbf{200}
			& \textbf{1.00}
			& 24.39
			& $4.07\times10^{-3}$ \\
			
			\bottomrule
		\end{tabular}
	}
	
\end{table}

\subsubsection*{Case 2: $\epsilon_1=\epsilon_2=10^{-5}$}
\begin{figure}[H]
	\centering
	\begin{minipage}{0.5\textwidth}
		\centering
		\includegraphics[width=\linewidth]{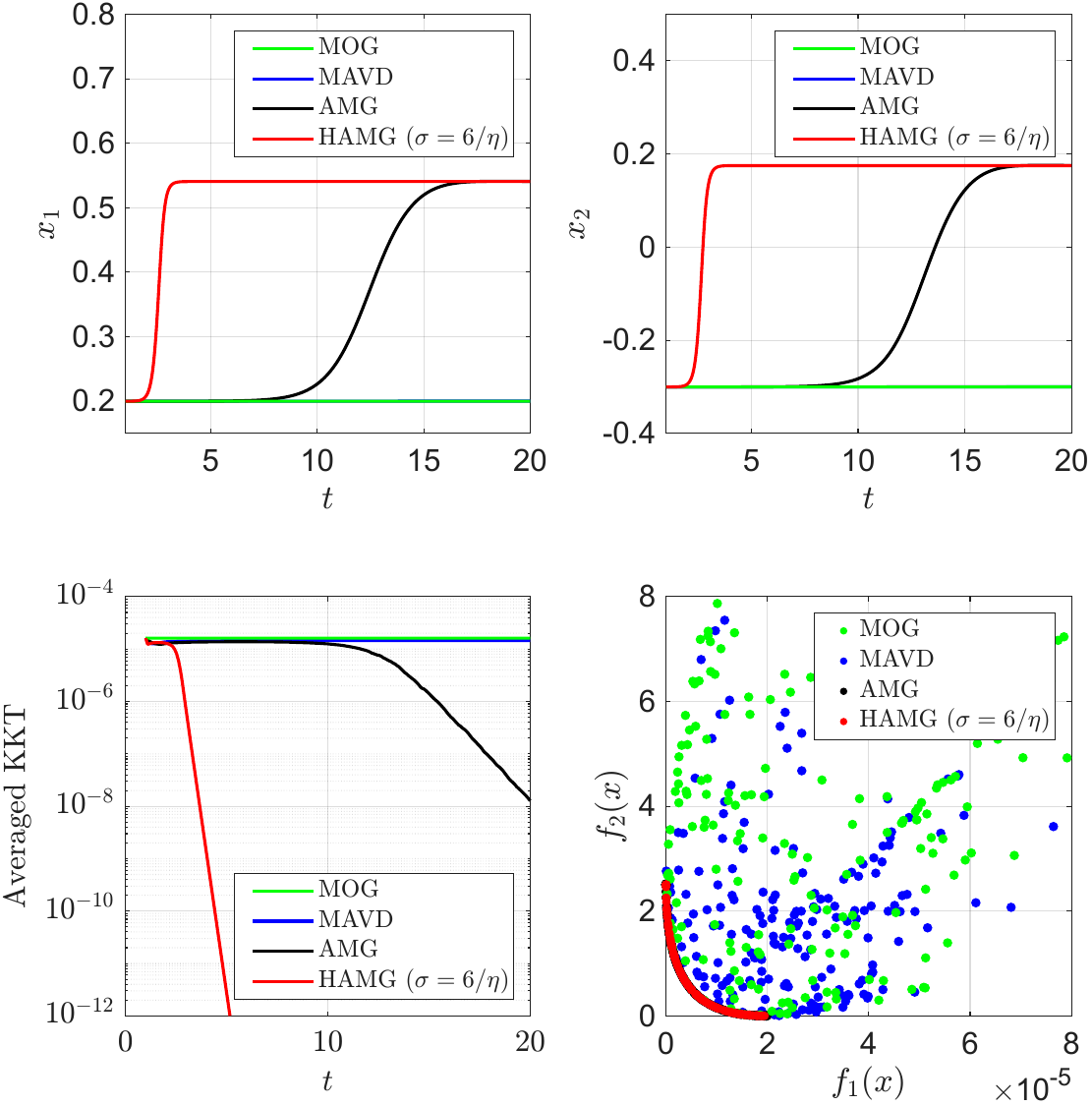}
	\end{minipage}
	
	\caption{Numerical results of \cref{eg:3} with $\epsilon_1=\epsilon_2=10^{-5}$.  The initial range is $[-2,2]^2$. Top line: the components of a sample trajectory. Bottom line: the averaged KKT residual (left) and  the approximate Pareto front (right).}
	\label{fig:img5}
\end{figure}

In this case, a common strong convexity parameter is still $\mu=10^{-5}$ and the two objectives have different scales. This makes the underlying multiobjective problem more difficult to solve. From \cref{tab:imbalanced example} and \cref{fig:img5}, we conclude that HAMG has shorter observed settling times, while AMG attains smaller terminal
IGD$^+$ values in both initial ranges. Their HV values agree at the displayed precision.
Both MOG and MAVD give poor performances on this imbalanced case.
\begin{table}[H]
	\centering
	\caption{Performances of different dynamical models for \cref{eg:3} with $\epsilon_1=\epsilon_2=10^{-5}$. ``$\infty$'' indicates that the corresponding residual does not satisfy the settling criterion by $T_{\max}=20$.}
	\label{tab:imbalanced example}
	\footnotesize
	\setlength{\tabcolsep}{10pt}
	\renewcommand{\arraystretch}{1.08}
	\resizebox{\textwidth}{!}{
		\begin{tabular}{lcccccc}
			\toprule
			Model
			& $T_{\mathrm{PF}}$
			& $T_{\mathrm{KKT}}$
			& NP
			& Purity $\uparrow$
			& HV $\uparrow$
			& IGD$^+$ $\downarrow$ \\
			\midrule
			
			\multicolumn{7}{c}{Initial range $[-2,2]^2$} \\
			\midrule
			
			MAVD
			& $\infty$
			& $\infty$
			& 16
			& 0.00
			& 3.75
			& $1.50\times10^{-4}$ \\
			
			AMG
			& $\infty$
			& 15.36
			& 180
			& 0.99
			& \textbf{3.76}
			& $\mathbf{5.33\times10^{-8}}$ \\
			
			MOG
			& $\infty$
			& $\infty$
			& 10
			& 0.10
			& 3.75
			& $1.05\times10^{-3}$ \\
			
			HAMG ($\sigma=6/\eta$)
			& \textbf{4.63}
			& \textbf{2.95}
			& \textbf{200}
			& \textbf{1.00}
			& \textbf{3.76}
			& $9.06\times10^{-7}$ \\
			
			\midrule
			\multicolumn{7}{c}{Initial range $[-50,50]^2$} \\
			\midrule
			
			MAVD
			& $\infty$
			& $\infty$
			& 2
			& 0.00
			& 53.68
			& $9.58\times10^{-1}$ \\
			
			AMG
			& $\infty$
			& 19.90
			& 86
			& 0.85
			& \textbf{55.94}
			& $\mathbf{7.48\times10^{-8}}$ \\
			
			MOG
			& $\infty$
			& $\infty$
			& 2
			& 0.00
			& 2.69
			& $5.18\times10^{1}$ \\
			
			HAMG ($\sigma=6/\eta$)
			& \textbf{5.50}
			& \textbf{3.51}
			& \textbf{200}
			& \textbf{1.00}
			& \textbf{55.94}
			& $1.45\times10^{-5}$ \\
			
			\bottomrule
		\end{tabular}
	}
\end{table}
\subsubsection{The log-sum-exp}
\begin{eg}\label{eg:convex}
	Consider the log-sum-exp objectives:	
	\begin{equation*} 
		f_j(x) =  
		\frac{\mu}{2}\|x\|^2
		+
		\ln
		\left(
		\sum_{i=1}^{p}
		\exp
		\left(
		\langle a_i^j,x\rangle-b_i^j
		\right)
		\right) 
		, \quad j=1,2,
	\end{equation*}	
	where $a^{j}_i \in \mathbb{R}^{n},\,b^{j}_i \in \mathbb{R}$ with $p=20,\,n=10$. 
\end{eg}	

\begin{figure}[H]
	\centering
	\begin{minipage}{0.5\textwidth}
		\centering
		\includegraphics[width=\linewidth]{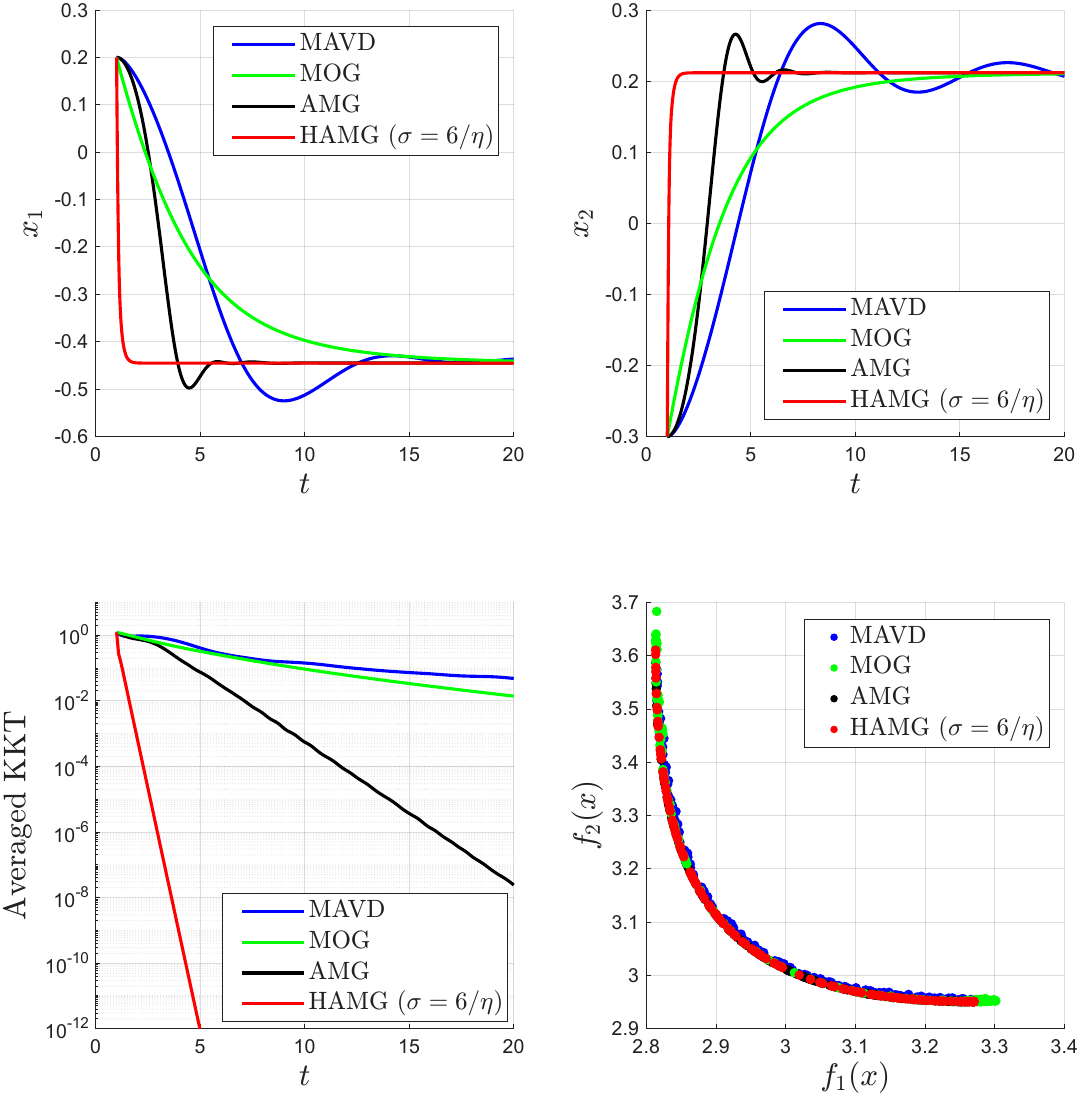}
	\end{minipage}
	
	\caption{Numerical results of \cref{eg:convex} with $\mu=0.1$.  The initial range is $[-2,2]^{10}$. Top line: the first component (left) and the second component (right) of a sample trajectory. Bottom line: the averaged KKT residual (left) and  the approximate Pareto front (right).}
	\label{fig:img6}
\end{figure}

\begin{table}[H]
	\centering
	\caption{Performances of different dynamical models for \cref{eg:convex} with $\mu=0.1$. ``$\infty$'' indicates that the corresponding residual does not
		satisfy the settling criterion by $T_{\max} = 20$.}
	\label{tab:performance_10d}
	\footnotesize
	\setlength{\tabcolsep}{10pt}
	\renewcommand{\arraystretch}{1.08}
	\resizebox{\textwidth}{!}{
		\begin{tabular}{lcccccc}
			\toprule
			Model
			& $T_{\mathrm{PF}}$
			& $T_{\mathrm{KKT}}$
			& NP
			& Purity $\uparrow$
			& HV $\uparrow$
			& IGD$^+$ $\downarrow$ \\
			\midrule
			
			\multicolumn{7}{c}{Initial range $[-2,2]^{10}$} \\
			\midrule
			
			MAVD
			& $\infty$
			& $\infty$
			& 97
			& 0.02
			& 2.39
			& $3.21\times10^{-3}$ \\
			
			AMG
			& 10.44
			& 16.26
			& \textbf{200}
			& \textbf{1.00}
			& \textbf{2.40}
			& $\mathbf{3.79\times10^{-4}}$ \\
			
			MOG
			& $\infty$
			& $\infty$
			& 95
			& 0.31
			& 2.39
			& $1.45\times10^{-3}$ \\
			
			HAMG ($\sigma=6/\eta$)
			& \textbf{2.68}
			& \textbf{3.01}
			& 189
			& \textbf{1.00}
			& \textbf{2.40}
			& $5.64\times10^{-4}$ \\
			
			\midrule
			\multicolumn{7}{c}{Initial range $[-50,50]^{10}$} \\
			\midrule
			
			MAVD
			& $\infty$
			& $\infty$
			& 2
			& 0.00
			& 5.04
			& $6.06\times10^{-1}$ \\
			
			AMG
			& 14.79
			& $\infty$
			& \textbf{200}
			& \textbf{1.00}
			& \textbf{7.70}
			& $\mathbf{4.03\times10^{-4}}$ \\
			
			MOG
			& $\infty$
			& $\infty$
			& 3
			& 0.00
			& 3.23
			& $1.33\times10^{0}$ \\
			
			HAMG ($\sigma=6/\eta$)
			& \textbf{3.49}
			& \textbf{3.55}
			& 169
			& \textbf{1.00}
			& \textbf{7.70}
			& $7.14\times10^{-4}$ \\
			
			\bottomrule
		\end{tabular}
	}
\end{table} 

\begin{figure}[H]
	\centering
	\begin{minipage}{0.5\textwidth}
		\centering
		\includegraphics[width=\linewidth]{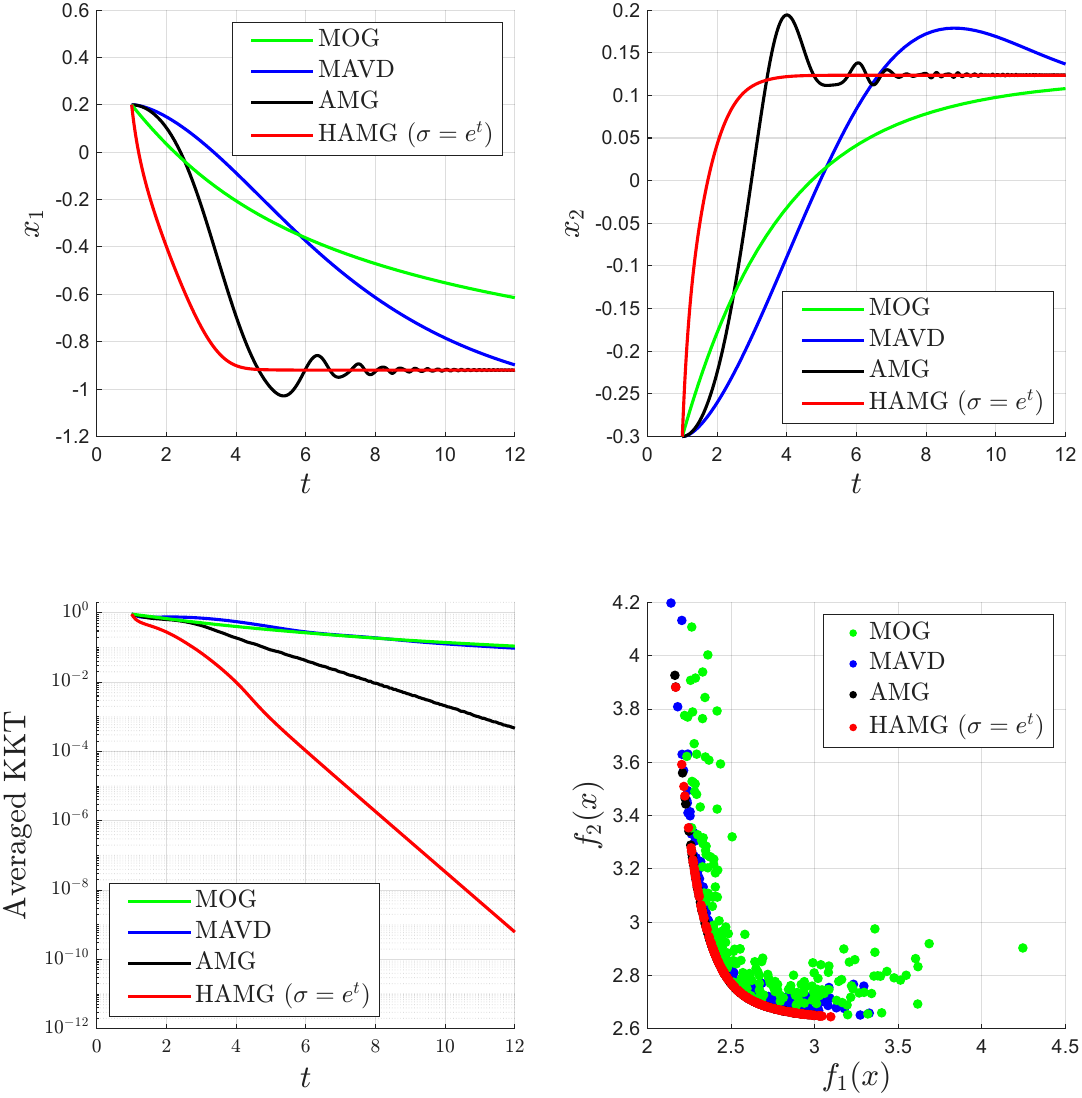}
	\end{minipage}
	
	\caption{Numerical results of \cref{eg:convex} with $\mu=0$.  The initial range is $[-2,2]^{10}$. Top line: the first component (left) and the second component (right) of a sample trajectory. Bottom line: the average  KKT residual (left) and  the approximate Pareto front (right).}
	\label{fig:convex}
\end{figure}

We consider two cases $\mu=0.1$ and $\mu=0$ separately, and the numerical results are reported in \cref{tab:performance_10d,fig:img6}, and \cref{tab:convex,fig:convex} respectively. As we can see, HAMG has smaller observed settling times for the residuals to reach the $\varepsilon$ level. In both cases, the Pareto approximation quantities of HAMG and AMG are more robust than others with respect to initial ranges. Especially, HAMG exhibits smoother trajectories and a visible reduction in oscillations, compared with the pronounced oscillations of AMG and MAVD.

\begin{table}[H]
	\centering
	\caption{Performance comparison of different models for \cref{eg:convex} with $\mu=0$. ``$\infty$'' indicates that the corresponding residual does not
		satisfy the settling criterion by $T_{\max} = 12$.}
	\label{tab:convex}
	\footnotesize
	\setlength{\tabcolsep}{10pt}
	\renewcommand{\arraystretch}{1.08}
	\resizebox{\textwidth}{!}{
		\begin{tabular}{lcccccc}
			\toprule
			Model
			& $T_{\mathrm{PF}}$
			& $T_{\mathrm{KKT}}$
			& NP
			& Purity $\uparrow$
			& HV $\uparrow$
			& IGD$^+$ $\downarrow$ \\
			\midrule
			
			\multicolumn{7}{c}{Initial range $[-2,2]^{10}$} \\
			\midrule
			
			MAVD
			& 12.00
			& $\infty$
			& 59
			& 0.03
			& 5.06
			& $1.37\times10^{-2}$ \\
			
			AMG
			& 11.33
			& $\infty$
			& \textbf{200}
			& \textbf{1.00}
			& 5.05
			& $\mathbf{1.00\times10^{-3}}$ \\
			
			MOG
			& 12.00
			& $\infty$
			& 28
			& 0.04
			& 5.04
			& $2.25\times10^{-2}$ \\
			
			HAMG ($\sigma=e^t$)
			& \textbf{5.10}
			& \textbf{8.17}
			& \textbf{200}
			& \textbf{1.00}
			& \textbf{5.07}
			& $\mathbf{1.00\times10^{-3}}$ \\
			
			\midrule
			\multicolumn{7}{c}{Initial range $[-50,50]^{10}$} \\
			\midrule
			
			MAVD
			& $\infty$
			& $\infty$
			& 3
			& 0.00
			& 1698.30
			& $2.37\times10^{1}$ \\
			
			AMG
			& 11.99
			& $\infty$
			& 145
			& 0.98
			& 3514.20
			& $\mathbf{1.24\times10^{-2}}$ \\
			
			MOG
			& $\infty$
			& $\infty$
			& 3
			& 0.00
			& 957.27
			& $4.51\times10^{1}$ \\
			
			HAMG ($\sigma=e^t$)
			& \textbf{7.02}
			& \textbf{9.83}
			& \textbf{174}
			& \textbf{1.00}
			& \textbf{3518.00}
			& $7.97\times10^{-2}$ \\
			
			\bottomrule
		\end{tabular}
	} 
\end{table} 
\subsubsection{A nonconvex  example}
\begin{eg}\label{eg:nonconvex}
	Consider the following two objectives
	\begin{equation*}
		\begin{aligned}
			f_1(x)
			&=
			\frac{1}{2}
			\left(
			\sqrt{1+(x_1+x_2)^2}
			+
			\sqrt{1+(x_1-x_2)^2}
			+
			x_1-x_2
			\right)
			+
			\lambda \exp\bigl(-(x_1-x_2)^2\bigr), \\[6pt]
			f_2(x)
			&=
			\frac{1}{2}
			\left(
			\sqrt{1+(x_1+x_2)^2}
			+
			\sqrt{1+(x_1-x_2)^2}
			-
			x_1+x_2
			\right)
			+
			\lambda \exp\bigl(-(x_1-x_2)^2\bigr),
		\end{aligned}
	\end{equation*}
	where $x=(x_1,x_2)\in\R^2$ and $\lambda=0.6$.
\end{eg}

\begin{figure}[H]
	\centering
	\begin{minipage}{0.5\textwidth}
		\centering
		\includegraphics[width=\linewidth]{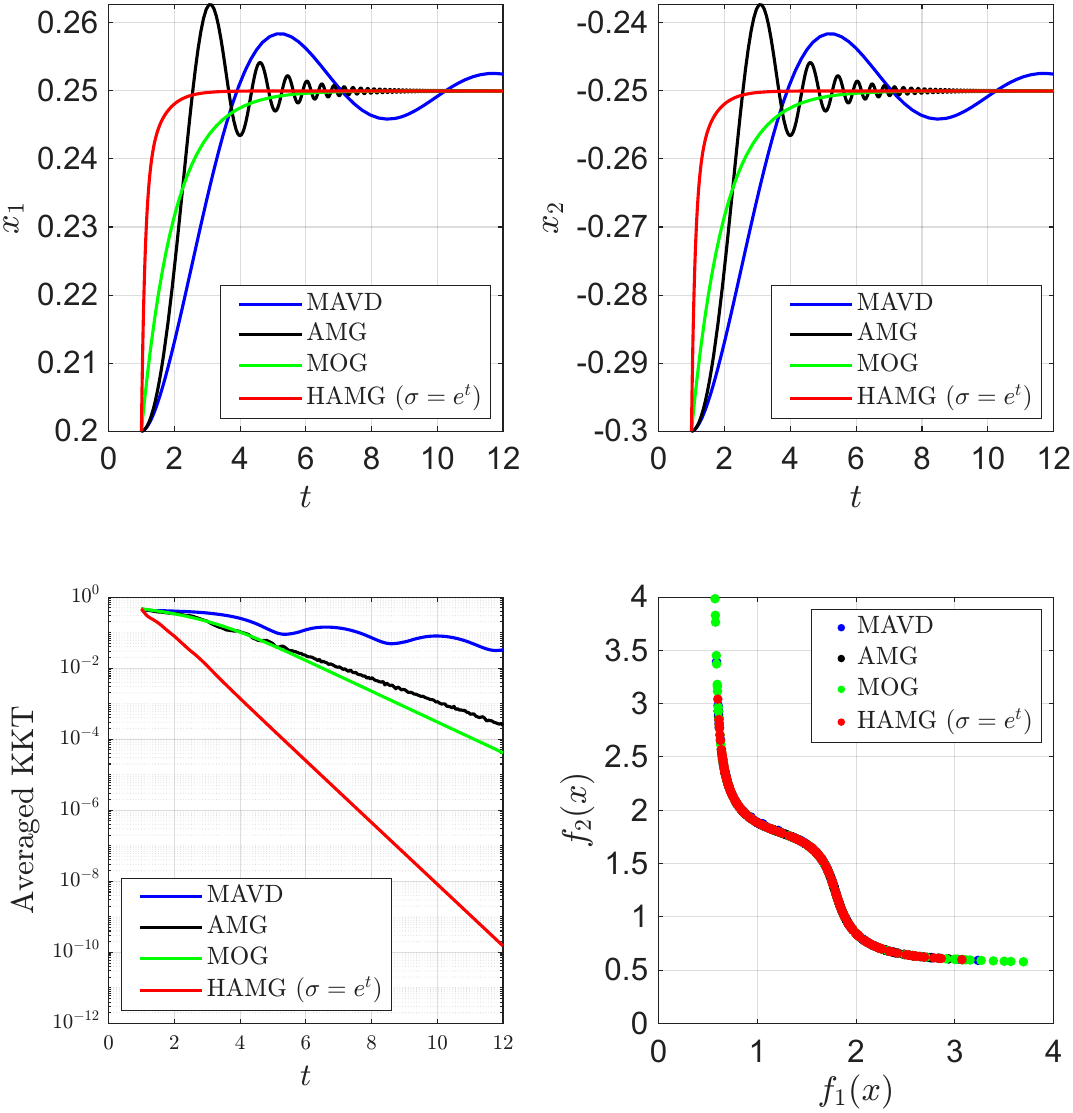}
	\end{minipage}
	\caption{Numerical results of \cref{eg:nonconvex}.
		  The initial range is $[-2,2]^2$. Top line: the components of the trajectory.
		Bottom line: the averaged KKT residual (left) and
		the approximate Pareto front (right).}
	\label{fig:img7}
\end{figure}
\begin{table}[H]
	\centering
	\caption{Performance comparison of different models for
		\cref{eg:nonconvex}. ``$\infty$'' indicates that the corresponding
		residual does not
		satisfy the settling criterion by $T_{\max} = 12$.}
	\label{tab:nonconvex}
	\footnotesize
	\setlength{\tabcolsep}{10pt}
	\renewcommand{\arraystretch}{1.08}
	\resizebox{\textwidth}{!}{
		\begin{tabular}{lcccccc}
			\toprule
			Model
			& $T_{\mathrm{PF}}$
			& $T_{\mathrm{KKT}}$
			& NP
			& Purity $\uparrow$
			& HV $\uparrow$
			& IGD$^+$ $\downarrow$ \\
			\midrule
			
			\multicolumn{7}{c}{Initial range $[-2,2]^2$} \\
			\midrule
			
			MAVD
			& $\infty$
			& $\infty$
			& 184
			& 0.80
			& 6.45
			& $5.60\times10^{-3}$ \\
			
			AMG
			& 9.30
			& $\infty$
			& \textbf{200}
			& \textbf{1.00}
			& 6.48
			& $4.06\times10^{-3}$ \\
			
			MOG
			& 7.54
			& $\infty$
			& \textbf{200}
			& \textbf{1.00}
			& \textbf{6.75}
			& $\mathbf{1.30\times10^{-3}}$ \\
			
			HAMG ($\sigma=e^t$)
			& \textbf{3.96}
			& \textbf{7.28}
			& \textbf{200}
			& \textbf{1.00}
			& 6.71
			& $1.34\times10^{-3}$ \\
			
			\midrule
			\multicolumn{7}{c}{Initial range $[-50,50]^2$} \\
			\midrule
			
			MAVD
			& $\infty$
			& $\infty$
			& 27
			& 0.15
			& 2263.30
			& $4.57\times10^{-2}$ \\
			
			AMG
			& $\infty$
			& $\infty$
			& 192
			& 0.97
			& \textbf{2264.60}
			& $3.42\times10^{-3}$ \\
			
			MOG
			& $\infty$
			& $\infty$
			& 66
			& 0.92
			& 2263.30
			& $5.54\times10^{-2}$ \\
			
			HAMG ($\sigma=e^t$)
			& \textbf{5.56}
			& \textbf{9.02}
			& \textbf{200}
			& \textbf{1.00}
			& \textbf{2264.60}
			& $\mathbf{1.70\times10^{-3}}$ \\
			
			\bottomrule
		\end{tabular}
	}
\end{table}
This experiment also lies outside the convexity assumptions of
our convergence theorems. Both objective functions are nonconvex,
and thus we set $\mu=0$ for HAMG and AMG.
The corresponding numerical results can be found in
\cref{tab:nonconvex} and \cref{fig:img7}.
For the small initial range, HAMG achieves the shortest
observed settling times, while the approximation quality of
the different methods is relatively close.
For the large initial range, HAMG outperforms all other
methods in terms of observed settling times and Pareto front
approximation quality within the prescribed maximal time.
In particular, HAMG achieves the highest purity and the smallest IGD$^+$ value, while its hypervolume ties with AMG at the displayed precision.


\section{Conclusion}\label{sec:conclu}

In this paper, we proposed an accelerated multiobjective gradient flow with implicit Hessian-driven
damping (HAMG). Its first-order formulation requires a single projection onto the convex hull
of the objective gradients. Proper time rescaling connects HAMG with the single-objective HNAG flow \cite{chen_first_2019} and yields multiobjective analogues of DIN \cite{alvarez2002second} and DIN-AVD \cite{attouch2016fast} for specific parameter choices. Under boundedness of the common lower level sets and
Lipschitz continuity of the gradients on bounded sets, we established a priori energy
estimate and a global existence result in finite dimensions. A Lyapunov analysis yields decay
estimates for a merit function and the KKT residual. With the additional restriction on the scaling parameters, every global solution converges weakly to a weak Pareto point, and the convergence is strong when all objectives are strongly convex. The numerical experiments illustrate rapid
residual decay and reduced oscillations on the tested problems.

%

An interesting future topic is the underlying discrete-time aspect, including multiobjective gradient algorithms and convergence analysis. The key is to find proper numerical schemes and parameter settings, together with efficient restart strategies, such that
the fast convergence property, especially the oscillation-reduction effect, can be retained after discretization.

%

\appendix
\section{Technical Lemmas}
\begin{lem}\label{lem:est-eta/A}
	Let $A(t)$ be defined by \cref{eq:U-def-second-order}, then $	\int_0^t
 \eta(\tau)/A^{2}(\tau)
	\dd\tau
	\leq
\eta_0/2$ for $t>0$.
\end{lem}
\begin{proof}
	Set $	G(t):=\frac{\eta(t)}{A^2(t)} $, then $G(0) = \eta_0$ and 
	\[
	G'(t) = \frac{A^2(t)\eta'(t)-2\eta(t)A(t)A'(t)}{A^4(t)} = \frac{A(t)\eta'(t)-2\eta(t)A'(t)}{A^3(t)}.
	\] 
	In view of \cref{eq:U-def-second-order}, we have \(A'(t)=(1+\mu\sigma(t))A(t)\) and
	by	\cref{eq:HAMG-eta},   
	\[
	\begin{aligned}
			A(t)\eta'(t)-2\eta(t)A'(t)={}&A(t)\left[\eta'(t)-2\eta(t)(1+\mu\sigma(t))\right]\\
			={}&-A(t)\left[\sigma(t)\eta^2(t)+\eta(t)(2+\mu\sigma(t))\right]\leq-2\eta(t)A(t),
	\end{aligned}
	\] 
	which implies immediately that $G'(t)\leq-2 G(t)$ and $G(t)\leq \eta_0 e^{-2t}$. Thus we obtain $	\int_0^t
	\eta(\tau)/A^2(\tau)
	\dd\tau
	\leq
	\eta_0
	\int_0^t e^{-2\tau}\dd\tau
	\leq
\eta_0/2$ for all $t>0$.
	This completes the proof.
\end{proof}

The following lemma follows from
\cite[Lemma~5.1]{AbbasAttouchSvaiter2014} and \cite[Lemma 5]{Bot2016}, with a simple integration argument.
\begin{lem}\label{lem:h-limit}
	Assume that (i) $g\in C[0,\infty) \cap L^1(0,\infty)$;  (ii) $a\in C[0,\infty)$ and $a(t)\geq0$ for all $t\geq0$.
	Let $h\in C^1[0,\infty)$ be nonnegative and satisfy $	h'(t)
	+  	a(t)h(t)
	\leq  g(t)$.
	Then $\lim_{t\to\infty}h(t)$ exists, and  if $\int_{0}^{\infty}a(t)\dd t=\infty$, then $\lim_{t\to\infty}h(t)
	= 0$.
\end{lem}
		
\bibliographystyle{abbrv}

\end{document}